\documentclass[seceqn,dvips]{arxaihp}
\usepackage{amssymb}
\usepackage{euscript}
\usepackage{graphicx}

\aid{0}
\pubyear{2008}
\volume{44}
\issue{1}
\firstpage{47}
\lastpage{88}
\doi{10.1214/07-AIHP106}

\makeatletter
\newcommand{\mathdsR}{\mathbb{R}}
\newcommand{\mathdsP}{\mathbb{P}}
\newcommand{\mathdsN}{\mathbb{N}}
\newcommand{\mathdsL}{\mathbb{L}}
\newcommand{\mathdsV}{\mathbb{V}}

\newcommand{\PP}{\mathrm{P}}
\def\<{\langle}
\def\>{\rangle}
\let\Sigma\varSigma
\let\Gamma\varGamma
\newtheorem{TEO}{Theorem}[section]

\newtheorem{lemma}[TEO]{Lemma}
\newtheorem{cor}[TEO]{Corollary}
\newproclaim{dfn}{Definition}[section]
\newproclaim{hypo}{Hypothesis}[section]
\newproclaim{cdtn}{Condition}[section]
\newproclaim{exm}{Example}[section]
\newproclaim{rmk}{Remark}[section]

\mathchardef\varTheta="0102
\mathchardef\varPi="0105
\mathchardef\varPsi="0109

\makeatother

\begin{document}
\begin{frontmatter}

\title{Iterative feature selection in least square regression estimation}
\runtitle{Iterative feature selection in least square regression estimation}

\begin{aug}
\author{\fnms{Pierre}\inits{P.} \snm{Alquier}\thanksref{}\ead[label=e1]{alquier@ensae.fr}}
\runauthor{P. Alquier}
\address{Laboratoire de Probabilit\'es et Mod\`{e}les Al\'{e}atoires,
Universit\'{e} Paris 6 and
Laboratoire de Statistique, Crest,
3, Avenue Pierre Larousse,
92240 Malakoff, France. \printead{e1}}
\end{aug}

\received{\sday{7} \smonth{1} \syear{2005}}
\revised{\sday{12} \smonth{12} \syear{2005}}
\accepted{\sday{23} \smonth{1} \syear{2006}}

%
\begin{abstract}
This paper presents a new algorithm to perform regression
estimation, in both the inductive and transductive setting. The
estimator is defined as a linear combination of functions in a given
dictionary. Coefficients of the combinations are computed
sequentially using projection on some simple sets. These sets are
defined as confidence regions provided by a deviation (PAC)
inequality on an estimator in one-dimensional models. We prove that
every projection the algorithm actually improves the performance of
the estimator. We give all the estimators and results at first in
the inductive case, where the algorithm requires the knowledge of
the distribution of the design, and then in the transductive case,
which seems a more natural application for this algorithm as we do
not need particular information on the distribution of the design in
this case. We finally show a connection with oracle inequalities,
making us able to prove that the estimator reaches minimax rates of
convergence in Sobolev and Besov spaces.
\end{abstract}
\begin{abstract}[language=french]
Cette article pr\'{e}sente un nouvel algorithme
d'estimation de r\'{e}gression, dans les contextes inductifs et
transductifs. L'estimateur est d\'{e}fini par une combinaison lin\'{e}aire
de fonctions choisies dans un dictionnaire donn\'{e}. Les coefficients
de cette combinaison sont calcul\'{e}s par des projections successives
sur des ensembles simples. Ces ensembles sont d\'{e}finis comme des
r\'{e}gions de confiance donn\'{e}es par une in\'{e}galit\'{e} de d\'
{e}viation (ou
in\'{e}galit\'{e} PAC). On d\'{e}montre en particulier que chaque
projection au
cours de l'algorithme am\'{e}liore effectivement l'estimateur obtenu. On
donne tout d'abord les r\'{e}sultats dans le contexte inductif, o\`{u}
l'algorithme n\'{e}cessite la connaissance de la distribution du design,
puis dans le contexte transductif, plus naturel ici puisque
l'algorithme s'applique sans la connaissance de cette distribution.
On \'{e}tablit finalement un lien avec les in\'{e}galit\'{e}s d'oracle,
permettant de montrer que notre estimateur atteint les vitesses
optimales dans les espaces de Sobolev et de Besov.
\end{abstract}

%
\begin{keyword}[class=MSC]
\kwd[Primary ]{62G08}
\kwd[; secondary ]{62G15}
\kwd{68T05}
\end{keyword}

\begin{keyword}
\kwd{Regression estimation}
\kwd{Statistical learning}
\kwd{Confidence regions}
\kwd{Thresholding methods}
\kwd{Support vector machines}
\end{keyword}

\end{frontmatter}

\section{The setting of the problem}

We give here notations and introduce the inductive and transductive settings.

\subsection{Transductive and inductive settings}

Let $(\mathcal{X},\mathcal{B})$ be a measure space and
let $\mathcal{B}_{\mathbb{R}}$ denote the Borel $\sigma$-algebra on
$\mathbb{R}$.

\subsubsection{The inductive setting}

In the inductive setting, we assume that $P$ is a distribution on pairs
$Z=(X,Y)$
taking values in $(\mathcal{X}\times\mathdsR ,\mathcal{B}\otimes
\mathcal{B}_{\mathdsR })$,
that $P$ is such that:
\[
P |Y | < \infty,
\]
and that we observe $N$ independent pairs $Z_{i}=(X_{i},Y_{i})$ for
$i\in\{1,\ldots,N\}$.
Our objective is then to estimate the regression function on the basis
of the observations.

\begin{dfn}[(The regression function)]
We denote:
\begin{eqnarray*}
\lefteqn{f\dvtx \mathcal{X}  \rightarrow\mathdsR ,}
\\
\lefteqn{x  \mapsto P(Y|X=x) .}
\end{eqnarray*}
\end{dfn}

\subsubsection{The transductive setting}

In the transductive case, we will assume that, for a given integer
$k>0$, $P_{(k+1)N}$ is some exchangeable probability measure on the
space $ ((\mathcal{X}\times\mathdsR )^{(k+1)N} ,
(\mathcal{B}\otimes\mathcal{B}_{\mathdsR })^{\otimes2N} )$.
We will write
$(X_{i},Y_{i})_{i=1,\ldots,(k+1)N}=(Z_{i})_{i=1,\ldots,(k+1)N}$ a random
vector distributed according to $P_{(k+1)N}$.

\begin{dfn}[(Exchangeable probability distribution)]
\label{exchp} For any integer $j$, let $\mathfrak{S}_{j}$ denote the
set of all permutations of $\{1,\ldots,j\}$. We say that $P_{(k+1)N}$
is exchangeable if for any $\sigma\in\mathfrak{S}_{(k+1)N}$ we have:
$(X_{\sigma(i)},Y_{\sigma(i)})_{i=1,\ldots,(k+1)N}$ has the same
distribution under $P_{(k+1)N}$ that
$(X_{i},Y_{i})_{i=1,\ldots,(k+1)N}$.
\end{dfn}

We assume that we observe $(X_{i},Y_{i})_{i=1,\ldots,N}$ and
$(X_{i})_{i=N+1,\ldots,(k+1)N}$; and the observation
$(X_{i},Y_{i})_{i=1,\ldots,(k+1)N}$ is usually called the training
sample, while the other part of the vector,
$(X_{i},Y_{i})_{i=N+1,\ldots,(k+1)N}$ is called the test sample. In
this case, we only focus on the estimation of the values
$(Y_{i})_{i=N+1,\ldots,(k+1)N}$. This is why Vapnik \cite{Vapnik}
called this kind of inference ``transductive inference'' when he
introduced it.

Note that in this setting, the pairs $(X_{i},Y_{i})$ are not
necessarily independent, but are identically distributed. We will
let $P$ denote their marginal distribution, and we can here again
define the regression function $f$.

Actually, most statistical problems being usually formulated in the
inductive setting, the reader may wonder about the pertinence of the
study of the transductive setting. Let us think of the following
examples: in quality control, or in a sample survey, we try to infer
informations about a whole population from observations on a small
sample. In this cases, transductive inference seems actually more
adapted than inductive inference, with $N$ the size of the sample
and $(k+1)N$ the size of the population. One can see that the use of
inductive results in this context is only motivated by the large
values of $k$ (the inductive case is the limit case of the
transductive case where $k\rightarrow+\infty$). In the problems
connected with regression estimation or classification, we can
imagine a case where a lot of images are collected for example on
the internet. The time to label every picture according to the fact
that it represents, or not, a given object being too long, one can
think of labeling only $1$ over $k+1$ images, and to use then a
transductive algorithm to label automatically the other data. We
hope that these examples can convince the reader that the use of the
transductive setting is not unrealistic. However, the reader that is
not convinced should remember that the transductive inference was
first introduced by Vapnik mainly as a tool to study the inductive
case: there are techniques to get rid of the second part of the
sample by taking an expectation with respect to it and obtain
results valid in the inductive setting (see for example a result by
Panchenko used in this paper, \cite{Panchenko}).

\subsection{The model}

In both settings, we are going to use the same model to estimate the
regression function: $\varTheta$. The only thing we assume about
$\varTheta$ is that it is a vector space of functions.

Note in particular that we do not assume that $f$ belongs to
$\varTheta$.

\subsection{Overview of the results}

In both settings, we give a PAC inequality on the risk of estimators
in one-dimensional models of the form:
\[
\bigl\{\alpha\theta(\cdot),\alpha\in\mathdsR  \bigr\}
\]
for a given $\theta\in\varTheta$.

This result motivates an algorithm that performs iterative feature
selection in order to perform regression estimation. We will then
remark that the selection procedure gives the guarantee that every
selected feature actually improves the current estimator.

In the inductive setting (Section~\ref{inductive}), it means that we
estimate $f(\cdot)$ by a function $\hat{\theta}(\cdot)\in\varTheta$, but the
selection procedure can only be performed if the statistician knows
the marginal distribution $P_{(X)}$ of $X$ under $P$.

In the transductive case (Section~\ref{sectiontrans}), the
estimation of $Y_{N+1},\ldots,Y_{(k+1)N}$ can be performed by the
procedure without any prior knowledge about the marginal
distribution of $X$ under $P$. We first focus on the case $k=1$, and
then on the general case $k\in\mathdsN ^{*}$.

Finally, in Section~\ref{rate}, we use the main result of the paper
(the fact that every selected feature improves the performance of
the estimator) as an oracle inequality, to compute the rate of
convergence of the estimator in Sobolev and Besov spaces.

The last section (Section~\ref{proofs}) is dedicated to the proofs.

The literature on iterative methods for regression estimation is
very important, let us mention one of the first algorithm, AdaLine,
by Widrow and Hoff \cite{Widrow}, or more recent versions like
boosting, see \cite{adaboost} and
the references within. The technique developed here has some
similarity with the so-called greedy algorithms, see
\cite{greedy} (and the references within) for a
survey and some recent results. However, note that in this
techniques, the iterative update of the estimator is motivated by
algorithmic issues, and is not motivated statistically. In
particular, AdaLine has no guarantee against overfitting if the
number of variables $m$ is large (say $m=N$). For greedy algorithms,
on has to specify a particular penalization if one wants to get a
guarantee against overfitting. The same remark can be done about
boosting algorithm. Here, the algorithm is motivated by a
statistical result, and as a consequence has theoretical guarantees
against overlearning. It stays however computationally feasible,
some pseudo-code is given in the paper.

Closer to our technique are the methods of aggregation of
statistical estimators, see \cite{Nemi} and \cite{Tsybagg}
and more recently the mirror descent algorithm studied
in \cite{Ju2} or \cite{Ju1}. In this papers, oracle
inequalities are given ensuring that the estimator performs as well
as the best (linear or convex) aggregation of functions in a given
family, up to an optimal term. Note that these inequalities are
given in expectation, here almost all results are given in a
deviation bound (or PAC bound, a bound that is true with high
probability, from which we derive a bound in expectation in Section~\ref{rate}). Similar bounds where given for the PAC-Bayesian model
aggregation developed by Catoni \cite{Cat7}, Yang \cite{Yang} and
Audibert~\cite{AudibertReg}. In some way, the algorithm proposed in
this paper can be seen as a practical way to implement these
results.

Note that nearly all the methods in the papers mentioned previously
where designed especially for the inductive setting. Very few
algorithms were created specifically for the transductive regression
problem. The algorithm described in this paper seems more adapted to
the transductive setting (remember that the procedure can be
performed in the inductive setting only if the statistician knows
the marginal distribution of $X$ under $P$, while there is no such
assumption in the transductive context).

Let us however start with a presentation of our method in the
inductive context.

\section{Main theorem in the inductive case, and application to estimation}
\label{inductive}

\subsection{Additional definition}

\begin{dfn}
We put:
\begin{eqnarray*}
\lefteqn{R(\theta)  = P \bigl[ \bigl(Y-\theta(X) \bigr)^{2} \bigr],}
\\
\lefteqn{r(\theta)  =\frac{1}{N}\sum_{i=1}^{N}
\bigl(Y_{i}-\theta(X_{i}) \bigr)^{2} ,}
\end{eqnarray*}
and in this setting, our objective is $\overline{\theta}$ given by:
\[
\overline{\theta} \in\mathop{\arg\min}_{\theta\in\varTheta} R(\theta) .
\]
\end{dfn}

\subsection{Main theorem}

We suppose that we have an integer $m\in\mathdsN $ and that we are
given a finite family of functions:
\[
\varTheta_{0} = \{\theta_{1},\ldots,\theta_{m}\} \subset\varTheta.
\]

\begin{dfn}
Let us put, for any $k \in\{1,\ldots,m\}$:
\begin{eqnarray*}
\lefteqn{\overline{\alpha}_{k}  = \mathop{\arg\min}_{\alpha\in\mathdsR }
R(\alpha\theta_{k})
= \frac{P [\theta_{k}(X)Y ]}{P [\theta_{k}(X)^{2} ]},}
\\
\lefteqn{\hat{\alpha}_{k}  = \mathop{\arg\min}_{\alpha\in\mathdsR }
r(\alpha\theta_{k})
= \frac{{(1/N)}\sum_{i=1}^{N} \theta_{k}(X_{i})Y_{i}}{(1/N)\sum_{i=1}^{N} \theta_{k}(X_{i})^{2}},}
\\
\lefteqn{\mathcal{C}_{k}  = \frac{(1/N)\sum_{i=1}^{N}
\theta_{k}(X_{i})^{2} }{P [\theta_{k}(X)^{2} ]} .}
\end{eqnarray*}
\end{dfn}

\begin{TEO} \label{lastTH}
Moreover, let us assume that $P$ is such that $|f|$ is bounded by a
constant $B$, and such that:
\[
P \bigl\{\bigl[Y-f(X)\bigr]^{2} \bigr\} \leq\sigma^{2} < +\infty.
\]
We have, for any $\varepsilon>0$, with $P^{\otimes N}$-probability
at least $1-\varepsilon$, for any $k\in\{1,\ldots,m\}$:
%
\begin{equation}
\label{lasteq}
R(\mathcal{C}_{k}\hat{\alpha}_{k}\theta_{k})
- R(\overline{\alpha}_{k}\theta_{k}) \leq
\frac{4 [1+\log(2m/\varepsilon) ]}{N} \biggl[\frac{
(1/N)\sum_{i=1}^{N}\theta_{k}(X_{i})^{2}Y_{i}^{2}
}{P[\theta_{k}(X)^{2}]} + B^{2} + \sigma^{2} \biggr].
\end{equation}
\end{TEO}

The proof of this theorem is given in Section~\ref{subprooflastTH}. 

\subsection{Application to regression estimation}
\label{sectionestimation}

\subsubsection[Interpretation of Theorem 2.1 in terms of
confidence intervals]{Interpretation of Theorem \textup{{\protect\ref{lastTH}}} in terms of
confidence intervals}

\begin{dfn}
Let us put, for any $ (\theta,\theta') \in\varTheta^{2}$:
\[
d_{P}\bigl(\theta,\theta'\bigr)
= \sqrt{P_{(X)} \bigl[ \bigl(\theta(X)-\theta'(X) \bigr)^{2} \bigr]}.
\]
Let also $\|\cdot\|_{P}$ denote the norm associated with this
distance, $ \|\theta\|_{P}=d_{P}(\theta,0)$, and
$ \langle\cdot,\cdot\rangle_{P}$ the associated scalar product:
\[
\bigl\langle\theta,\theta' \bigr\rangle _{P}=P \bigl[ \theta(X)\theta'(X) \bigr] .
\]
\end{dfn}

Because $\overline{\alpha}_{k}=\mathop{\arg\min}_{\alpha\in\mathdsR }
R(\alpha\theta_{k})$ we have:
\[
R(\mathcal{C}_{k}\hat{\alpha}_{k}\theta_{k})-R(\overline{\alpha
}_{k}\theta_{k})
=
d_{P}^{2}(\mathcal{C}_{k}\hat{\alpha}_{k}\theta_{k},\overline{\alpha
}_{k}\theta_{k})
.
\]

So the theorem can be written:
\[
P^{\otimes N} \bigl\{\forall k\in\{1,\ldots,m\},
d_{P}^{2}(\mathcal{C}_{k}\hat{\alpha}_{k}\theta_{k},\overline{\alpha
}_{k}\theta_{k})
\leq\beta(\varepsilon,k) \bigr\} \geq1-\varepsilon,
\]
where $\beta(\varepsilon,k)$ is the right-hand side of inequality~(\ref{lasteq}).

Now, note that $\overline{\alpha}_{k}\theta_{k}$ is the orthogonal
projection of:
\[
\overline{\theta} = \mathop{\arg\min}_{\theta\in\varTheta} R(\theta)
\]
onto the space $\{\alpha\theta_{k},\alpha\in\mathdsR \}$, with
respect to the inner product $ \langle\cdot,\cdot\rangle_{P}$:
\[
\overline{\alpha}_{k} = \mathop{\arg\min}_{\alpha\in\mathdsR } d_{P} (\alpha
\theta_{k},\overline{\theta} ).
\]

\begin{dfn}
We define, for any $k$ and $\varepsilon$:
\[
\mathcal{CR}(k,\varepsilon)= \biggl\{\theta\in\varTheta\dvt
\bigg| \bigg\<\theta-\mathcal{C}_{k}\hat{\alpha}_{k}\theta_{k},\frac{\theta_{k}}{
\|\theta_{k} \|_{P}} \biggr\>_{P} \bigg|
\leq\sqrt{\beta(\varepsilon,k)} \biggr\} .
\]
\end{dfn}

Then the theorem is equivalent to the following corollary.

\begin{cor}
We have:
\[
P^{\otimes N} \bigl[\forall k\in\{1,\ldots,m\}, \overline{\theta} \in\mathcal
{CR}(k,\varepsilon) \bigr]
\geq1-\varepsilon.
\]
\end{cor}

In other words: $\bigcap_{k\in\{1,\ldots,m\}}
\mathcal{CR}(k,\varepsilon)$ is a confidence region at level
$\varepsilon$ for $ \overline{\theta}$.

\begin{dfn}
We write $\varPi^{k,\varepsilon}_{P}$ the orthogonal projection into
$\mathcal{CR}(k,\varepsilon)$ with respect to the distance $d_{P}$.
\end{dfn}

Note that this orthogonal projection is not a projection on a linear
subspace of $\varTheta$, and so it is not a linear mapping.

\subsubsection{The algorithm}

The previous corollaries of Theorem \ref{lastTH} motivate the
following iterative algorithm:
\begin{itemize}
\item choose $\theta^{(0)}\in\varTheta$, for example, $\theta^{(0)}=0$;

\item at step $n\in\mathdsN ^{*}$, we have: $\theta^{(0)},\ldots,\theta
^{(n-1)}$.
Choose $k(n)\in\{1,\ldots,m\}$ (this choice can of course be data
dependent), and take:
\[
\theta^{(n)}=\varPi_{P}^{k(n),\varepsilon}\theta^{(n-1)} ;
\]

\item we can use the following stopping rule: $ \|\theta^{(n-1)}-\theta
^{(n)} \|_{P}^{2}\leq\kappa$, where $0<\kappa<\frac{1}{N}$.
\end{itemize}

\begin{dfn} \label{hatf}
Let $n_{0}$ denote the stopping step, and:
\[
\hat{\theta}(\cdot) = \theta^{(n_{0})}(\cdot)
\]
the corresponding function.
\end{dfn}

\subsubsection{Results and comments on the algorithm}

\begin{TEO} \label{propalgo}
We have:
\[
P^{\otimes N} \bigl[\forall n\in\{1,\ldots,n_{0}\},
R \bigl(\theta^{ (n )} \bigr) \leq
R \bigl(\theta^{ (n-1 )} \bigr) -
d_{P}^{2} \bigl(\theta^{(n)},\theta^{(n-1)} \bigr)
\bigr]
\geq1-\varepsilon.
\]
\end{TEO}

\begin{pf}
This is just a consequence of the preceding corollary. Let us assume
that:
\[
 \forall k\in\{1,\ldots,m\},\quad
R(\mathcal{C}_{k}\hat{\alpha}_{k}\theta_{k})
- R(\overline{\alpha}_{k}\theta_{k}) \leq\beta(\varepsilon,k).
\]
Let us choose $n\in\{1,\ldots,n_{0}\}$. We have, for a
$k\in\{1,\ldots,m\}$:
\[
\theta^{(n)}=\varPi_{P}^{k,\varepsilon}\theta^{(n-1)},
\]
where $\varPi_{P}^{k,\varepsilon}$ is the projection into a convex set
that contains $\overline{\theta}$. This implies that:
\[
\bigl\<\theta^{(n)}-\theta^{(n-1)},\overline{\theta}-\theta^{(n)} \bigr\>_{P} \geq
0 ,
\]
or:
\[
d_{P}^{2} \bigl(\theta^{(n-1)},\overline{\theta} \,\bigr) \geq d_{P}^{2} \bigl(\theta
^{(n)},\overline{\theta} \,\bigr)
+ d_{P}^{2} \bigl(\theta^{(n-1)},\theta^{(n)} \bigr) ,
\]
which can be written:
\[
R \bigl[\theta^{(n-1)} \bigr]-R(\overline{\theta}) \geq R \bigl[\theta^{(n)}
\bigr]-R(\overline{\theta})
+ d_{P}^{2} \bigl(\theta^{(n-1)},\theta^{(n)} \bigr).
\]\upqed
\end{pf}

Actually, the main point in the motivation of the algorithm is that,
with probability at least $1-\varepsilon$, whatever the current
value $\theta^{(n)}\in\varTheta$, whatever the feature
$k\in\{1,\ldots,m\}$ (even chosen on the basis of the data),
$\varPi^{k,\varepsilon}_{P}\theta^{(n)}$ is a better estimator than
$\theta^{(n)}$.

So we can choose $k(n)$ as we want in the algorithm. For example,
Theorem \ref{propalgo} motivates the choice:
\[
k(n) = \mathop{\arg\max}_{k} d_{P}^{2} \bigl(\theta^{(n-1)},\mathcal{CR}(k,\varepsilon
) \bigr).
\]
This version of the algorithm is detailed in Fig.~1. If looking
for the exact maximum of
\[
d_{P} \bigl(\theta^{(n-1)},\mathcal{CR}(k,\varepsilon) \bigr)
\]
with respect to $k$ is too computationally intensive we can use any
heuristic to choose $k(n)$, or even skip this maximization and take:
\[
k(1) = 1,\ldots,\qquad k(m) = m,\qquad k(m+1)= 1,\ldots,\qquad k(2m) =
m,\ldots.
\]

\begin{figure}
\begin{tabular}{p{11.0cm}}
\hline
We have $\varepsilon>0$, $\kappa>0$, $N$ observations
$(X_{1},Y_{1}),\ldots,(X_{N},Y_{N})$, $m$ features
$\theta_{1}(\cdot),\ldots,\theta_{m}(\cdot)$ and
$c=(c_{1},\ldots,c_{m})=(0,\ldots,0)\in\mathdsR ^{m}$. Compute at first
every $\hat{\alpha}_{k}$ and $\beta(\varepsilon,k)$ for
$k\in\{1,\ldots,m\}$. Set $n \leftarrow0$.
\bigskip

Repeat:
\begin{itemize}
\item set $n\leftarrow n+1$;

\item set \textit{best}\_\textit{improvement}${}\leftarrow0$;

\item for $k\in\{1,\ldots,m\}$, compute:
\begin{eqnarray*}
\lefteqn{v_{k} = P \bigl[\theta_{k}(X)^{2} \bigr] ,}
\\
\lefteqn{\gamma_{k} \leftarrow\hat{\alpha}_{k} -
\frac{1}{v_{k}}\sum_{j=1}^{m}c_{j}P \bigl[\theta_{j}(X)\theta_{k}(X)
\bigr],}
\\
\lefteqn{\delta_{k} \leftarrow v_{k} \bigl( |\gamma_{k} | - \beta(\varepsilon,k)
\bigr)_{+}^{2},}
\end{eqnarray*}
and if $\delta_{k}>\mbox{\textit{best}\_\textit{improvemen}}t$, set:
\begin{eqnarray*}
\lefteqn{\mbox{best\_improvement} \leftarrow\delta_{k},}
\\
\lefteqn{k(n) \leftarrow k ;}
\end{eqnarray*}

\item if \textit{best}\_\textit{improvement}${}>0$ set:
\[
c_{k(n)} \leftarrow c_{k(n)} + \operatorname{sgn}(\gamma_{k(n)}) \bigl( |\gamma_{k(n)} | -
\beta\bigl(\varepsilon,k(n)\bigr) \bigr)_{+} ;
\]
\end{itemize}
until \textit{best}\_\textit{improvement}${}<\kappa$ (where $\operatorname{sgn}(x)=-1$ if $x\leq0$ and
$1$ otherwise).

\bigskip

Note that at each step $n$, $\theta^{(n)}$ is given by:
\[
\theta^{(n)} (\cdot) = \sum_{k=1}^{m} c_{k} \theta_{k}(\cdot) ,
\]
so after the last step we can return the estimator:
\[
\hat{\theta}(\cdot) = \sum_{k=1}^{m} c_{k} \theta_{k} (\cdot) .
\]
\\
\hline
\end{tabular}
\caption{Detailed version of the feature selection algorithm.}
\end{figure}

\begin{exm}
Let us assume that $\mathcal{X}= [0,1]$ and let us put
$\varTheta=\mathdsL_{2}(P_{(X)})$. Let
$(\theta_{k})_{k\in\mathdsN ^{*}}$ be an orthonormal basis
of~$\varTheta$. The choice of $m$ should not be a problem, the algorithm
itself avoiding itself overlearning we can take a large value of $m$
like $m=N$. In this setting, the algorithm is a procedure for (soft)
thresholding of coefficients. In the particular case of a wavelets
basis, see \cite{Donoho} or \cite{Wavelets} for a presentation of wavelets coefficient
thresholding. Here, the threshold is not necessarily the same for
every coefficient. We can remark that the sequential projection on
every $k$ is sufficient here:
\[
k(1) = 1,\ldots,k(m) = m,
\]
after that $\theta^{(m+n)}=\theta^{(m)}$ for every $n\in\mathdsN $
(because all the directions of the different projections are
orthogonals).
\end{exm}

Actually, it is possible to prove that the estimator is able to
adapt itself to the regularity of the function to achieve a good
mean rate of convergence. More precisely, if we assume that the true
regression function has an (unknown) regularity $\beta$, then it is
possible to choose $m$ and $\varepsilon$ in such a way that the rate
of convergence is:
\[
N^{{-2\beta}/{(2\beta+1)}}\log N .
\]
We prove this point in Section~\ref{rate}.

\begin{rmk}
Note that in its general form, the algorithm does not require any
assumption about the dictionary of functions
$\varTheta_{0}=\{\theta_{1},\ldots,\theta_{m}\}$. This family can be
non-orthogonal, it can even be redundant (the dimension of the
vector space generated by $\varTheta_{0}$ can be smaller than $m$).
\end{rmk}

\begin{rmk}
It is possible to generalize Theorem \ref{lastTH} to models of
dimension larger than $1$. The algorithm itself can take advantage
of these generalizations. This point is developed in
\cite{AlqThese}, where some experiences about the performances of
our algorithm can also be found.
\end{rmk}

\subsection[Additional notations for some refinements of Theorem 2.1]{Additional notations for some refinements of Theorem \textup{{\protect\ref{lastTH}}}}

Note that an improvement of the inequality in Theorem \ref{lastTH}
(inequality~(\ref{lasteq})) would allow to apply
the same method, but would lead to smaller confidence regions and so
to better performances. The end of this section is dedicated to
improvements (and generalizations) of this bound.

\begin{hypo*}
Until the end of Section~\ref{inductive}, we
assume that $\varTheta$ and $P$ are such that:
\[
\forall\theta\in\varTheta,\quad   P\exp\bigl[\theta(X)Y \bigr] < +\infty.
\]
\end{hypo*}

\begin{dfn}
For any random variable $T$ we put:
\begin{eqnarray*}
\lefteqn{V(T)  = P \bigl[ (T-PT )^{2} \bigr],}
\\
\lefteqn{M^{3}(T)  = P \bigl[ (T-PT )^{3} \bigr],}
\end{eqnarray*}
and we define, for any $\gamma\geq0$, $P_{ \gamma T}$ by:
\[
\frac{ \mathrm{d}P_{ \gamma T}}{\mathrm{d}P} = \frac{\exp (\gamma T ) }{ P [\exp
(\gamma T ) ]} .
\]
For any random variables $T,T'$ and any $\gamma\geq0$ we put:
\begin{eqnarray*}
\lefteqn{V_{ \gamma T}\bigl(T'\bigr)  = P_{ \gamma T} \bigl[ \bigl(T'-P_{ \gamma T}T' \bigr)^{2} \bigr],}
\\
\lefteqn{M^{3}_{ \gamma T}\bigl(T'\bigr)  = P_{ \gamma T} \bigl[ \bigl(T'-P_{ \gamma T}T' \bigr)^{3}
\bigr].}
\end{eqnarray*}
\end{dfn}

Section~\ref{sub1} gives an improvement of Theorem \ref{lastTH}
while Section~\ref{subsvm} extends it to the case of a
data-dependant family $\varTheta_{0}$.

\subsection[Refinements of Theorem 2.1]{Refinements of Theorem \textup{{\protect\ref{lastTH}}}}
\label{sub1}

\begin{TEO} \label{TH1}
Let us put:
\[
W_{\theta} = \theta(X)Y-P \bigl(\theta(X)Y \bigr) .
\]
Then we have, for any $\varepsilon>0$, with $P^{\otimes N}$-probability
at least $1-\varepsilon$,
for any $k\in\{1,\ldots,m\}$:
\[
R(\mathcal{C}_{k}\hat{\alpha}_{k}\theta_{k})
- R(\,\overline{\alpha}_{k}\theta_{k}) \leq
\frac{2\log (2m/\varepsilon)}{N}
\frac{V (W_{\theta_{k}} )}{P [\theta_{k}(X)^{2} ]}
+ \frac{\log^{3} (2m/\varepsilon)}{N^{{3}/{2}}}
C_{N}(P,m,\varepsilon,\theta_{k}) ,
\]
where we have:
\begin{eqnarray*}
C_{N}(P,m,\varepsilon,\theta_{k}) &=&
I_{\theta_{k}} \biggl(\sqrt{\frac{2\log ({2m}/{\varepsilon})}{N
V (W_{\theta_{k}} )}} \biggr)^{2} \frac{\sqrt{2}
}{V (W_{\theta_{k}} )^{{5}/{2}}P [\theta_{k}(X)^{2} ]}
\\
&&{} + I_{\theta_{k}} \biggl(\sqrt{\frac{2\log ({2m}/{\varepsilon})}{N
V (W_{\theta_{k}} )}} \biggr)^{4}
\frac{\log^{2} ({2m}/{\varepsilon})
}{\sqrt{N}V (W_{\theta_{k}} )^{6}P [\theta_{k}(X)^{2} ]},
\end{eqnarray*}
with:
\[
I_{\theta} (\gamma ) = \int_{0}^{1}(1-\beta)^{2}M^{3}_{\beta\gamma
W_{\theta}} (W_{\theta} )\,\mathrm{d}\beta.
\]
\end{TEO}

For the proof, see Section~\ref{subproof1}.

Actually, the method we proposed requires to be able to compute
explicitly the upper bound in this theorem. Remark that, with
$\varepsilon$ and $m$ fixed:
\[
C_{N}(P,m,\varepsilon,\theta_{k}) \mathop{\longrightarrow}_{N \rightarrow+\infty}
\frac{\sqrt{2} [M^{3} (W_{\theta_{k}} ) ]^{2} }{9
V (W_{\theta_{k}} )^{{5}/{2}}P [\theta_{k}(X)^{2} ]}
,
\]
and so we can choose to consider only the first-order term. Another
possible choice is to make stronger assumptions on $P$ and
$\varTheta_{0}$ that allow to upper bound explicitly $
C_{N}(P,m,\varepsilon,\theta_{k})$. For example, if we assume that
$Y$ is bounded by $C_{Y}$ and that $\theta_{k}(\cdot)$ is bounded by
$C_{k}'$ then $W_{\theta_{k}}$ is bounded by $C_{k}=2C_{Y}C_{k}'$
and we have (basically):
\[
C_{N}(P,m,\varepsilon,\theta_{k}) \leq\frac{ 64\sqrt{2} C_{k}^{2} }{9
V (W_{\theta_{k}} )^{{5}/{2}}P [\theta_{k}(X)^{2} ]}
+ \frac{4096 C_{k}^{4} \log^{3} ({2m}/{\varepsilon})}{81
\sqrt{N}V (W_{\theta_{k}} )^{6}
P [\theta_{k}(X)^{2} ]} .
\]

The main problem is actually that the first-order term contains the quantity
$V (W_{\theta_{k}} )$ that is not observable, and we would like to be able
to replace this quantity by its natural estimator:
\[
\hat{V}_{k} = \frac{1}{N}\sum_{i=1}^{N} \Biggl[Y_{i} \theta_{k}(X_{i}) - \frac
{1}{N}\sum_{j=1}^{N} Y_{j} \theta_{k}(X_{j}) \Biggr]^{2} .
\]

The following theorem justifies this method.

\begin{TEO} \label{TH1bis}
If we assume that there is a constant $c$ such that:
\[
\forall k\in\{1,\ldots,m\},\quad P \bigl[\exp\bigl(c W_{\theta_{k}}^{2} \bigr) \bigr]<\infty,
\]
we have, for any $\varepsilon>0$, with $P^{\otimes N}$-probability at
least $1-\varepsilon$,
for any $k\in\{1,\ldots,m\}$:
\[
R(\mathcal{C}_{k}\hat{\alpha}_{k}\theta_{k})
- R(\overline{\alpha}_{k}\theta_{k}) \leq
\frac{2\log ({4m}/{\varepsilon})}{N}
\frac{\hat{V}_{k}}{P [\theta_{k}(X)^{2} ]} +
\frac{\log ({4m}/{\varepsilon})}{N^{{3}/{2}}}
C_{N}'(P,m,\varepsilon,\theta_{k}) ,
\]
where we have:
\[
\hat{V}_{k} = \frac{1}{N}\sum_{i=1}^{N} \Biggl[Y_{i}\theta_{k}(X_{i}) - \frac
{1}{N}\sum_{j=1}^{N} Y_{j}\theta_{k}(X_{j}) \Biggr]^{2},
\]
and
\begin{eqnarray*}
C_{N}'(P,m,\varepsilon,\theta_{k}) &=&  C_{N} \biggl(P,m,\frac{\varepsilon
}{2},\theta_{k} \biggr) \log^{2} \frac{4m}{\varepsilon}
\\
&&{}+
\frac{2\log^{{1}/{2}}(2m/\varepsilon)}{P [\theta_{k}(X)^{2} ]}
\biggl[\sqrt{2V \bigl(W_{\theta_{k}}^{2} \bigr)} +
\frac{\log (2m/\varepsilon)}{\sqrt{N}
V (W_{\theta_{k}}^{2} ) }
J_{\theta_{k}} \biggl(\sqrt{\frac{2\log(2m/\varepsilon)}{N
V (W_{\theta_{k}}^{2} )}} \biggr) \biggr]
\\
&&{}+
\frac{2\log^{{1}/{2}}({4m}/{\varepsilon})}{P [\theta_{k}(X)^{2} ]}
\biggl[ \sqrt{2V (W_{\theta_{k}} )} +
\frac{\log^{2}({2m}/{\varepsilon})}{\sqrt{N}
V (W_{\theta_{k}} )^{3} }
I_{\theta_{k}} \biggl(\sqrt{\frac{2\log({4m}/{\varepsilon})}{N
V (W_{\theta_{k}} )}} \biggr) \biggr]
\\
&&{}\times \Biggl[ \Bigg|\frac{2}{N}\sum_{i=1}^{N}Y_{i} \theta_{k}(X_{i}) \Bigg|
\sqrt{\frac{2V (W_{\theta_{k}} )\log (4m/\varepsilon)}{N}}
+ \frac{\log^{{5}/{2}} ({2m}/{\varepsilon})}{N
V (W_{\theta_{k}} )^{3} }
I_{\theta_{k}} \biggl(\sqrt{\frac{2\log ({4m}/{\varepsilon})}{N
V (W_{\theta_{k}} )}} \biggr) \Biggr]
\end{eqnarray*}
and
\[
J_{\theta} (\gamma) = \int_{0}^{1}(1-\beta)^{2}M^{3}_{\gamma\beta
W_{\theta_{k}}^{2}} \bigl(W_{\theta}^{2} \bigr)\,\mathrm{d}\beta.
\]
\end{TEO}

The proof is given in Section~\ref{subproof1}.

\subsection{An extension to the case of Support Vector Machines}

\label{subsvm}

Thanks to a method due to Seeger \cite{Seeger}, it is possible to
extend this method to the case where
the set $\varTheta_{0}$ is data dependent in the following way:
\[
\varTheta_{0}(Z_{1},\ldots,Z_{N},N) = \bigcup_{i=1}^{N} \varTheta
_{0}(Z_{i},N) ,
\]
where for any $z\in\mathcal{X}\times\mathdsR $, the cardinality of
the set $\varTheta_{0}(z,N)$ depends only on $N$, not on $z$. We will
write $m'(N)$ this cardinality. So we have:
\[
\big|\varTheta_{0}(Z_{1},\ldots,Z_{N},N) \big|\leq N \big|\varTheta_{0}(Z_{i},N) \big| = N
m'(N) .
\]
We put:
\[
\varTheta_{0}(Z_{i},N) = \{\theta_{i,1},\ldots,\theta_{i,m'(N)} \} .
\]

In this case, we need some adaptations of our previous notations.

\begin{dfn}
We put, for $i\in\{1,\ldots,N\}$:
\[
r_{i}(\theta) = \frac{1}{N-1} \mathop{\sum_{j\in\{1,\ldots,N\},}}_{j \neq
i}  \bigl(Y_{j}- \theta(X_{j}) \bigr)^{2} .
\]
For any $(i,k) \in\{1,\ldots,N\}\times\{1,\ldots,m'(N)\}$, we write:
\begin{eqnarray*}
\lefteqn{\hat{\alpha}_{i,k}  = \mathop{\arg\min}_{\alpha\in\mathdsR } r_{i}(\alpha\theta_{i,k})
= \frac{ \sum_{j\neq i
} \theta_{i,k}(X_{j})Y_{j} }{ \sum_{ j\neq i }
\theta_{i,k}(X_{j})^{2}
},} \\
\lefteqn{\overline{\alpha}_{i,k}  = \mathop{\arg\min}_{\alpha\in\mathdsR } R(\alpha
\theta_{i,k})
= \frac{P [\theta_{i,k}(X)Y ]}{P [\theta_{i,k}(X)^{2} ]},}
\\
\lefteqn{\mathcal{C}_{i,k}  = \frac{{1}/{(N-1)}\sum_{ j\neq i }
\theta_{i,k}(X_{j})^{2} }{P [\theta_{i,k}(X)^{2} ]} .}
\end{eqnarray*}
\end{dfn}

\begin{TEO} \label{THSVM}
We have, for any $\varepsilon>0$, with $P^{\otimes N}$-probability at
least $1-\varepsilon$,
for any $k\in\{1,\ldots,m'(N)\}$ and $i\in\{1,\ldots,N\}$:
\begin{eqnarray*}
R(\mathcal{C}_{i,k}\hat{\alpha}_{i,k}\theta_{i,k})
- R(\overline{\alpha}_{i,k}\theta_{i,k}) &\leq&
\frac{2\log ({2Nm'(N)}/{\varepsilon})}{N-1}
\frac{V (W_{\theta_{i,k}} )}{P [\theta_{i,k}(X)^{2} ]}
\\
&&{} + \frac{\log^{3} ({2Nm'(N)}/{\varepsilon})}{(N-1)^{{3}/{2}}}
C_{N-1}\bigl(P,Nm'(N),\varepsilon,\theta_{i,k}\bigr).
\end{eqnarray*}
\end{TEO}

The proof is given in Section~\ref{subproof1}.

We can use this theorem to build an estimator using the algorithm
described in the previous subsection,
with obvious changes in the notations.

\begin{exm}
Let us consider the case where $\mathcal{H}$ is a Hilbert space with
scalar product $ \langle\cdot,\cdot\rangle$, and:
\[
\varTheta= \bigl\{\theta(\cdot) = \bigl\<h,\varPsi(\cdot) \bigr\>,  h\in\mathcal{H}
\bigr\}
\]
where $\varPsi$ is an application
$\mathcal{X}\rightarrow\varTheta$. Let us put
$\varTheta_{0}[(x,y),N]= \{ \<\varPsi(x),\varPsi(\cdot) \> \}$.
In this case we have $m'(N)=1$ and the estimator is of the from:
\[
\hat{\theta}(\cdot) = \sum_{i=1}^{N} \alpha_{i,1} \bigl\<\varPsi(X_{i}),\varPsi(\cdot) \bigr\> .
\]
Let us define,
\[
K\bigl(x,x'\bigr)= \bigl\<\varPsi(x),\varPsi\bigl(x'\bigr) \bigr\> ,
\]
the function $K$ is called the kernel, and:
\[
I=\{1\leq i \leq N\dvt \alpha_{i,1} \neq0 \} ,
\]
that is called the set of support vectors.
Then the estimate has the form of a support vector machine (SVM):
\[
\hat{\theta}(\cdot) = \sum_{i\in I} \alpha_{i,1} K(X_{i},\cdot) .
\]
SVM where first introduced by Boser, Guyon and Vapnik \cite{SVM_FIRST}
in the context of classification,
and then generalized by Vapnik \cite{Vapnik} to the context of
regression estimation.
For a general introduction to SVM, see also \cite{Classif} and \cite{Cristianini}.
\end{exm}

\begin{exm} \label{MSVM}
A widely used kernel is the Gaussian kernel:
\[
K_{\gamma}\bigl(x,x'\bigr)= \exp\biggl(-\gamma\frac{d^{2}(x,x')}{2} \biggr) ,
\]
where $d(\cdot,\cdot)$ is some distance over the space $\mathcal{X}$ and
$\gamma>0$. But in practice, the choice of the parameter $\gamma$ is
difficult. A way to solve this problem is to introduce multiscale
SVM. We simply take $\varTheta$ as the set of all bounded functions
$\mathcal{X}\rightarrow\mathdsR $.
Now, let us put:
\[
\varTheta_{0}\bigl[(x,y),N\bigr]= \bigl\{K_{2}(x,\cdot),K_{2^{2}}(X,\cdot),\ldots
,K_{2^{m'(N)}}(x,\cdot) \bigr\}.
\]
In this case, we obtain an estimator of the form:
\[
\hat{\theta}(\cdot) = \sum_{k=1}^{m'(N)} \sum_{i\in I_{k}} \alpha_{i,k}
K_{2^{k}}(X_{i},\cdot) ,
\]
that could be called multiscale SVM. Remark that we can use this
technique to define SVM using simultaneously different kernels (not
necessarily the same kernel at different scales).
\end{exm}

\section{The transductive case}
\label{sectiontrans}

\subsection{Notations}

Let us recall that we assume that $k\in\mathdsN ^{*}$, that
$P_{(k+1)N}$ is some exchangeable probability measure (let us recall
that exchangeability is defined in Definition \ref{exchp}) on the space
$ ((\mathcal{X}\times\mathdsR )^{(k+1)N} ,
(\mathcal{B}\times\mathcal{B}_{\mathdsR })^{\otimes
(k+1)N} )$. Let
$(X_{i},Y_{i})_{i=1,\ldots,(k+1)N}=(Z_{i})_{i=1,\ldots,(k+1)N}$ denote a
random vector distributed according to $P_{(k+1)N}$.

Let us remark that under this condition, the marginal distribution
of every $Z_{i}$ is the same, we will call $P$ this distribution. In
the particular case where the observations are i.i.d., we will have
$P_{(k+1)N}=P^{\otimes(k+1)N}$, but what follows still holds for
general exchangeable distributions $P_{(k+1)N}$.

We assume that we observe $(X_{i},Y_{i})_{i=1,\ldots,N}$ and
$(X_{i})_{i=N+1,\ldots,(k+1)N}$. In this case, we only focus on the
estimation of the values $(Y_{i})_{i=N+1,\ldots,(k+1)N}$.

\begin{dfn}
We put, for any $\theta\in\varTheta$:
\begin{eqnarray*}
\lefteqn{r_{1}(\theta)  = \frac{1}{N}\sum_{i=1}^{N} \bigl(Y_{i}-\theta(X_{i})
\bigr)^{2},}
\\
\lefteqn{r_{2}(\theta)  =
\frac{1}{kN}\sum_{i=N+1}^{(k+1)N} \bigl(Y_{i}-\theta(X_{i})
\bigr)^{2}.}
\end{eqnarray*}
Our objective is:
\[
\overline{\theta}_{2} = \mathop{\arg\min}_{\theta\in\varTheta} r_{2}(\theta) ,
\]
if the minimum of $r_{2}$ is not unique then we take for $\overline
{\theta}_{2}$
any element of $\varTheta$ reaching the minimum value of $r_{2}$.
\end{dfn}

Let $\varTheta_{0}$ be a finite family of vectors belonging to
$\varTheta$, so that $|\varTheta_{0}|=m$. Actually, $\varTheta_{0}$ is
allowed to be data-dependent:
\[
\varTheta_{0} = \varTheta_{0}(X_{1},\ldots,X_{(k+1)N}),
\]
but we assume that the function $(x_{1},\ldots,x_{(k+1)N}) \mapsto
\varTheta_{0}(x_{1},\ldots,x_{(k+1)N})$ is exchangeable with respect to
its $(k+1)N$ arguments, and is such that $m=m(N)$ depends only on
$N$, not on $(X_{1},\ldots,X_{(k+1)N})$.

The problem of the indexation of the elements of $\varTheta_{0}$ is not
straightforward and we must be very careful about it.
Let $<_{\varTheta}$ be a complete order on $\varTheta$,
and write:
\[
\varTheta_{0} = \{\theta_{1},\ldots,\theta_{m} \},
\]
where
\[
\theta_{1}<_{\varTheta}\cdots<_{\varTheta}\theta_{m}.
\]
Remark that, in this case, every $\theta_{h}$ is an exchangeable
function of $(X_{1},\ldots,X_{(k+1)N})$.

\begin{dfn}
Now, let us write, for any $h \in\{1,\ldots,m\}$:
\begin{eqnarray*}
\lefteqn{\alpha_{1}^{h}  = \mathop{\arg\min}_{\alpha\in\mathdsR }
r_{1}(\alpha\theta_{h})
= \frac{\sum_{i=1}^{N} \theta_{h}(X_{i})Y_{i}}{\sum_{i=1}^{N} \theta
_{h}(X_{i})^{2}},}
\\
\lefteqn{\alpha_{2}^{h}  = \mathop{\arg\min}_{\alpha\in\mathdsR }
r_{2}(\alpha\theta_{h})
= \frac{\sum_{i=N+1}^{(k+1)N} \theta_{h}(X_{i})Y_{i}}{\sum
_{i=N+1}^{(k+1)N} \theta_{h}(X_{i})^{2}},}
\\
\lefteqn{\mathcal{C}^{h}  = \frac{(1/N)\sum_{i=1}^{N}
\theta_{h}(X_{i})^{2} }{(1/(kN))\sum_{i=N+1}^{(k+1)N}
\theta_{h}(X_{i})^{2} } .}
\end{eqnarray*}
\end{dfn}

\subsection{Basic results for $k=1$}

In a first time we focus on the case where $k=1$ as a method due to
Catoni \cite{Classif} brings a substantial simplification of the
bound in this case.

\begin{TEO} \label{thtrans}
We have, for any $ \varepsilon>0$, with $P_{2N}$-probability at
least $1-\varepsilon$, for any $h\in\{1,\ldots,m\}$:
\[
r_{2} \bigl[\bigl(\mathcal{C}^{h}
\alpha_{1}^{h}\bigr)\cdot \theta_{h} \bigr]-r_{2}\bigl(\alpha_{2}^{h}\cdot \theta_{h}\bigr)
\leq4 \biggl[
\frac{(1/N)\sum_{i=1}^{2N}\theta_{h}(X_{i})^{2}Y_{i}^{2}}
{(1/N)\sum_{i=N+1}^{2N}\theta_{h}(X_{i})^{2}} \biggr]
\frac{\log ({2m}/{\varepsilon})}{N} .
\]
\end{TEO}

\begin{rmk}
Here again, it is possible to make some hypothesis in order to make the
right-hand side of the theorem observable.
In particular, if we assume that:
\[
\exists B \in\mathdsR _{+},\quad   P\bigl(|Y| \leq B\bigr) =1 ,
\]
then we can get a looser observable upper bound:
\begin{eqnarray*}
\hspace*{-12pt}
P_{2N} \biggl\{ \forall k\in\{1,\ldots,m\},
r_{2} \bigl[\bigl(\mathcal{C}^{h}
\alpha_{1}^{h}\bigr)\cdot\theta_{h} \bigr]-r_{2}\bigl(\alpha_{2}^{h}\cdot\theta_{h}\bigr)
\leq4 \biggl[ B^{2}+
\frac{({1}/{N})\sum_{i=1}^{N}\theta_{h}(X_{i})^{2}Y_{i}^{2}}
{({1}/{N})\sum_{i=N+1}^{2N}\theta_{h}(X_{i})^{2}} \biggr]
\frac{\log ({2m}/{\varepsilon})}{N} \biggr\} \geq1-\varepsilon.
\end{eqnarray*}
If we do not want to make this assumption, we can use the following
variant, that gives a first-order approximation for the bound.
\end{rmk}

\begin{TEO} \label{cortrans}
For any $ \varepsilon>0$, with $P_{2N}$-probability at least
$1-\varepsilon$, for any $h\in\{1,\ldots,m\}$:
\begin{eqnarray*}
\lefteqn{r_{2} \bigl[\bigl(\mathcal{C}^{h}
\alpha_{1}^{h}\bigr)\cdot\theta_{h}
\bigr]-r_{2}\bigl(\alpha_{2}^{h}\cdot\theta_{h}\bigr)}
\\
\lefteqn{\quad \leq\frac{8 \log ({4m}/{\varepsilon})}{N} \Biggl[
\frac{({1}/{N})\sum_{i=1}^{N}\theta_{h}(X_{i})^{2}Y_{i}^{2}}
{({1}/{N})\sum_{i=N+1}^{2N}\theta_{h}(X_{i})^{2}} +
\sqrt{\frac{({1}/{N})\sum_{i=1}^{2N}\theta_{h}(X_{i})^{4}Y_{i}^{4}\log
({2m}/{\varepsilon})}{2N}}\Biggr].}
\end{eqnarray*}
\end{TEO}

\begin{rmk}
Let us assume that $Y$ is such that we know two constants $b_{Y}$ and
$B_{Y}$ such that:
\[
P\exp\bigl(b_{Y} |Y| \bigr)\leq B_{Y} < \infty.
\]
Then we have, with probability at least $1-\varepsilon$:
\[
\sup_{i\in\{1,\ldots,2N\}} |Y_{i}| \leq\frac{1}{b_{Y}} \log\frac
{2NB_{Y}}{\varepsilon} .
\]
Combining both inequalities leads by a union bound argument leads
to:
\begin{eqnarray*}
\lefteqn{r_{2} \bigl[\bigl(\mathcal{C}^{h}
\alpha_{1}^{h}\bigr)\cdot\theta_{h}
\bigr]-r_{2}\bigl(\alpha_{2}^{h}\cdot\theta_{h}\bigr)}
\\
\lefteqn{\quad\leq\frac{8 \log ({8m}/{\varepsilon})}{N} \Biggl[
\frac{({1}/{N})\sum_{i=1}^{N}\theta_{h}(X_{i})^{2}Y_{i}^{2}}
{({1}/{N})\sum_{i=N+1}^{2N}\theta_{h}(X_{i})^{2}} +
\sqrt{\frac{({1}/{N})\sum_{i=1}^{2N}\theta_{h}(X_{i})^{4}\log ({4m}/{\varepsilon})\log^{4}
({4NB_{Y}}/{\varepsilon})}{2Nb_{Y}^{4}}} \Biggr].}
\end{eqnarray*}
\end{rmk}

The proofs of both theorems are given in the proofs section, more
precisely in Section~\ref{subproof2}.

Let us compare the first-order term of this theorem to the analogous
term in the inductive case (Theorems \ref{TH1}~and~\ref{TH1bis}).
The factor of the variance term is $8$ instead of $2$ in the
inductive case. A factor $2$ is to be lost because we have here the
variance of a sample of size $2N$ instead of $N$ in the inductive
case. But another factor $2$ is lost here. Moreover, in the
inductive case, we obtained the real variance of $Y\theta_{h}(X)$
instead of the moment of order $2$ here.

In the next subsection, we give several improvements of these
bounds, that allows to recover a real variance, and to recover the
factor $2$. We also give a version that allows to deal with a test
sample of different size, this being a generalization of Theorem
\ref{thtrans} more than of its improved variants.

We then give the analog of the algorithm proposed in the inductive
case in this transductive setting.

\subsection{Improvements of the bound and general values for $k$}

The proof of all the theorems of this subsection is given in the next section.

\subsubsection{Variance term (in the case $k=1$)}

We introduce some new notations.

\begin{dfn}
We write:
\[
\forall\theta\in\varTheta, r_{1,2}(\theta)=r_{1}(\theta)+r_{2}(\theta)
\]
and, in the case of a model $k\in\{1,\ldots,m\}$:
\[
\alpha_{1,2}^{h}=\mathop{\arg\min}_{\alpha\in\mathdsR } r_{1,2}(\alpha\theta
_{h}) .
\]
\end{dfn}

The we have the following theorem.

\begin{TEO}\label{thimp1}
We have, for any $ \varepsilon>0$, with $P_{2N}$-probability at
least $1-\varepsilon$, for any $h\in\{1,\ldots,m\}$:
\[
r_{2}\bigl(\mathcal{C}^{h}\alpha_{1}^{h}\theta_{h}\bigr)-r_{2}\bigl(\alpha
_{2}^{h}\theta_{h}\bigr)
\leq
4 \biggl[\frac{({1}/{N})\sum_{i=1}^{2N} [\theta_{h}(X_{i})Y_{i}-\alpha
_{1,2}^{h}\theta_{h}(X_{i})^{2} ]^{2}}
{({1}/{N})\sum_{i=N+1}^{2N}\theta_{h}(X_{i})^{2}} \biggr]\frac{\log({2m}/{\varepsilon})}{N}.
\]
\end{TEO}

For the proof see Section~\ref{subproof3}.

It is moreover possible to modify the upper bound to make it observable.
We obtain that with $P_{2N}$-probability at least $1-\varepsilon$, for any
$k\in\{1,\ldots,m\}$:
\begin{eqnarray*}
r_{2} \bigl[\bigl(\mathcal{C}^{h}
\alpha_{1}^{h}\bigr)\theta_{h} \bigr]-r_{2}\bigl(\alpha_{2}^{h}\theta_{h}\bigr)
\leq\frac{16 \log ({4m}/{\varepsilon})}{N}
\Biggl[\frac{1}{N}\sum_{i=1}^{N} \bigl(\theta_{h}(X_{i})Y_{i}-\alpha_{1}^{h}\theta
_{h}(X_{i})^{2} \bigr)^{2} \Biggr]
+ \EuScript{O} \biggl( \biggl[\frac{\log ({m}/{\varepsilon})}{N} \biggr]^{{3}/{2}} \biggr).
\end{eqnarray*}

So we can see that this theorem is an improvement on Theorem
\ref{thtrans} when some features $\theta_{h}(X)$ are well correlated
with $Y$. But we loose another factor $2$ by making the first-order
term of the bound observable.

\subsubsection{Improvement of the variance term ($k=1$)}

\begin{TEO} \label{thimp2}
We have, for any $ \varepsilon>0$, with $P_{2N}$-probability at
least $1-\varepsilon$, for any $h\in\{1,\ldots,m\}$:
\[
r_{2}\bigl(\mathcal{C}^{h}\alpha_{1}^{h}\theta_{h}\bigr)-r_{2}\bigl(\alpha
_{2}^{h}\theta_{h}\bigr)
\leq
\biggl[\frac{1}{1-{2\log ({2m}/{\varepsilon})}/{N}} \biggr]\frac{2\log
({2m}/{\varepsilon})}{N}
\frac{V_{1}(\theta_{h})+V_{2}(\theta_{h})}{({1}/{N})\sum_{i=N+1}^{2N}
\theta_{h}(X_{i})^{2}} ,
\]
where:
\begin{eqnarray*}
\lefteqn{V_{1}(\theta_{h})  = \frac{1}{N}\sum_{i=1}^{N} \Biggl[Y_{i}\theta
_{h}(X_{i})-\frac{1}{N}\sum_{j=1}^{N}Y_{j}\theta_{h}(X_{j})
\Biggr]^{2},}
\\
\lefteqn{V_{2}(\theta_{h})  =
\frac{1}{N}\sum_{i=N+1}^{2N} \Biggl[Y_{i}\theta_{h}(X_{i})-\frac{1}{N}\sum
_{j=N+1}^{2N}Y_{j}\theta_{h}(X_{j}) \Biggr]^{2}.}
\end{eqnarray*}
It is moreover possible to give an observable upper bound: we obtain
that with $P_{2N}$-probability at least $1-\varepsilon$, for any
$k\in\{1,\ldots,m\}$:
\begin{eqnarray*}
r_{2} \bigl[\bigl(\mathcal{C}^{h}
\alpha_{1}^{h}\bigr)\theta_{h} \bigr]-r_{2}\bigl(\alpha_{2}^{h}\theta_{h}\bigr)
&\leq& \biggl[\frac{1}{1-{2\log ({4m}/{\varepsilon})}/{N}} \biggr]
\frac{4\log ({4m}/{\varepsilon})}{N}\frac{V_{1}(\theta_{h})}{({1}/{N})\sum_{i=N+1}^{2N}
\theta_{h}(X_{i})^{2}}
\\
&&{} + \biggl[\frac{1}{1- {2\log ({4m}/{\varepsilon})}/{N}} \biggr]
2 \bigl(2+\sqrt{2} \bigr) \biggl(\frac{\log ({6m}/{\varepsilon})}{N} \biggr)^{{3}/{2}}
\frac{\sqrt{({1}/{N})\sum_{i=1}^{2N}\theta
_{h}(X_{i})^{4}Y_{i}^{4}}}{({1}/{N})\sum_{i=N+1}^{2N}
\theta_{h}(X_{i})^{2}} .
\end{eqnarray*}
\end{TEO}

Here again, we can make the bound fully observable under an
exponential moment or boundedness assumption about $Y$. For a
complete proof see Section~\ref{subproof4}.

\subsubsection{The general case ($k\in\mathdsN ^{*}$)}

We need some new notations in this case.

\begin{dfn}
Let us put:
\[
\mathbf{P} = \frac{1}{(k+1)N}\sum_{i=1}^{(k+1)N}\delta_{Z_{i}},
\]
and, for any $\theta\in\varTheta$:
\[
\mathdsV _{\theta} = \mathbf{P} \bigl\{ \bigl[ \bigl(\theta(X)Y \bigr)-\mathbf{P} \bigl(\theta
(X)Y \bigr) \bigr]^{2} \bigr\} .
\]
\end{dfn}

Then we have the following theorem.

\begin{TEO} \label{thimp3}
Let us assume that we have constants $B_{h}$ and $\beta_{h}$ such that,
for any $h\in\{1,\ldots,m\}$:
\[
P \exp\bigl(\beta_{h} \big|\theta_{h}(X_{i})Y_{i} \big| \bigr) \leq B_{h} .
\]
For any $\varepsilon>0$, with $P_{(k+1)N}$ probability at least
$1-\varepsilon$ we have,
for any $h\in\{1,\ldots,m\}$:
\begin{eqnarray*}
\lefteqn{r_{2} \bigl(\mathcal{C}^{h}\alpha_{1}^{h}\theta_{h} \bigr) -
r_{2} \bigl(\alpha_{2}^{h}\theta_{h} \bigr)}
\\
\lefteqn{\quad \leq
\frac{ (1+1/{k} )^{2}}{(1/(kN))\sum_{i=N+1}^{(k+1)N}
\theta_{h}(X_{i})^{2}} \biggl[ \frac{2 \mathdsV _{\theta_{h}}
\log ({4m}/{\varepsilon})}{N}}
\\
\lefteqn{\qquad{} + \frac{16 (\log ({4m}/{\varepsilon}) )^{{3}/{2}} (\log
(4(k+1)mNB_{h}/\varepsilon)
)^{3}}{3 \beta_{h}^{3} N^{{3}/{2}}
\mathdsV _{\theta_{h}}^{{1}/{2}} }
+ \frac{64 (\log(4m/\varepsilon) )^{2} (\log (4(k+1)mNB_{h}/{\varepsilon}))^{6}}{ 9 \beta_{h}^{6} N^{2}
\mathdsV _{\theta_{h}}^{2} }\biggr].}
\end{eqnarray*}
\end{TEO}

Here again, it is possible to replace the variance term by its natural
estimator:
\[
\hat{\mathdsV }_{\theta_{h}} =\frac{1}{N}\sum_{i=1}^{N} \Biggl[\theta
_{h}(X_{i})Y_{i}-\frac{1}{N}\sum_{j=1}^{N}\theta_{h}(X_{j})Y_{j} \Biggr]^{2} .
\]
For a complete proof of the theorem see the section dedicated to the
proofs (more precisely Section~\ref{subproof5}).

\subsection{Application to transductive regression}

We give here the interpretation of the preceding theorems in terms of
confidence; this
motivates an algorithm similar to the one described in the inductive case.

\begin{dfn}
We take, for any
$(\theta,\theta')\in\varTheta^{2}$:
\[
d_{2}\bigl(\theta,\theta'\bigr)
= \sqrt{\frac{1}{kN}\sum_{i=N+1}^{(k+1)N}
\bigl[\theta(X_{i})-\theta'(X_{i}) \bigr]^{2}}.
\]
Let also
$\|\theta\|_{2}=d_{2}(\theta,0)$ and:
\[
\bigl\<\theta,\theta' \bigr\>_{2}=\frac{1}{(k+1)N}\sum_{i=N+1}^{(k+1)N} \theta
(X_{i})\theta'(X_{i}) .
\]
We define, for any $h\in\{1,\ldots,m\}$ and $\varepsilon$:
\[
\mathcal{CR}(h,\varepsilon)= \bigl\{\theta\in\varTheta\dvt  \big|\bigl\<\theta-\mathcal
{C}^{h}\alpha_{1}^{h}\theta_{h},\theta_{h} \bigr\>_{2} \big|
\leq\sqrt{\beta(\varepsilon,h)} \bigr\} ,
\]
where $\beta(\varepsilon,h)$ is the upper bound in Theorem \ref{thtrans}
(or in any other theorem given in the transductive section).
\end{dfn}

For the same reasons as in the inductive case, these theorems imply
the following result.

\begin{cor} We have:
\[
P_{2N} \bigl[\forall h\in\{1,\ldots,m\}, \overline{\theta}_{2}\in\mathcal
{CR}(h,\varepsilon) \bigr]
\geq1-\varepsilon.
\]
\end{cor}

\begin{dfn}
We call $\varPi_{2}^{h,\varepsilon}$ the orthogonal projection into
$\mathcal{CR}(h,\varepsilon)$ with respect to the distance
$d_{2}$.
\end{dfn}

We propose the following algorithm:
\begin{itemize}
\item choose $\theta^{(0)}\in\varTheta$ (for example $0$);

\item at step $n\in\mathdsN ^{*}$, we have: $\theta^{(0)},\ldots,\theta
^{(n-1)}$. Choose $h(n)$, for example:
\[
h(n)=\mathop{\arg\max}_{h\in\{1,\ldots,m\}}d_{2}\bigl(\theta^{(n-1)},\mathcal
{CR}(h,\varepsilon)\bigr),
\]
and take:
\[
\theta^{(n)}=\varPi_{2}^{h(n),\varepsilon}\theta^{(n-1)} ;
\]

\item we can use the following stopping rule: $ \|\theta^{(n-1)}-\theta
^{(n)} \|_{2}^{2}\leq\kappa$
where $0<\kappa<\frac{1}{N}$.
\end{itemize}

\begin{dfn}
We write $n_{0}$ the stopping step, and:
\[
\theta(\cdot)= \theta^{(n_{0})}(\cdot)
\]
the corresponding function.
\end{dfn}

Here again we give a detailed version of the algorithm, see Fig.~2.
Remark that as in the inductive case, we are allowed to use
whatever heuristic to choose $k(n)$ if we want to avoid the
maximization.

\begin{figure}
\begin{tabular}{p{11.0cm}}
\hline
We have $\varepsilon>0$, $\kappa>0$, $N$ observations
$(X_{1},Y_{1}),\ldots,(X_{N},Y_{N})$ and also $X_{N+1},\ldots,X_{(k+1)N}$,
$m$ features $\theta_{1}(\cdot),\ldots,\theta_{m}(\cdot)$ and
$c=(c_{1},\ldots,c_{m})=(0,\ldots,0)\in\mathdsR ^{m}$. First, compute
every $\alpha_{1}^{h}$ and $\beta(\varepsilon,h)$ for
$h\in\{1,\ldots,m\}$. Set $n \leftarrow0$.

\bigskip

Repeat:
\begin{itemize}
\item set $n\leftarrow n+1$;

\item set \textit{best}\_\textit{improvement}${}\leftarrow0$;

\item for $h\in\{1,\ldots,m\}$, compute:
\begin{eqnarray*}
\lefteqn{v_{h} =
\frac{1}{kN}\sum_{i=N+1}^{(k+1)N}\theta_{h}(X_{i})^{2},}
\\
\lefteqn{\gamma_{h} \leftarrow\alpha^{h}_{1} -
\frac{1}{v_{h}}
\sum_{j=1}^{m}c_{j}\frac{1}{kN}\sum_{i=N+1}^{(k+1)N}\theta
_{j}(X_{i})\theta_{h}(X_{i}),}
\\
\lefteqn{\delta_{h} \leftarrow v_{h} \bigl( |\gamma_{h} | - \beta(\varepsilon,h)
\bigr)_{+}^{2},}
\end{eqnarray*}
and if $\delta_{h}>{}$\textit{best}\_\textit{improvement}, set:
\begin{eqnarray*}
\mbox{\textit{best\_improvement}} \leftarrow\delta_{h},
\\
h(n) \leftarrow h ;
\end{eqnarray*}

\item if \textit{best}\_\textit{improvement}${}>0$ set:
\[
c_{h(n)} \leftarrow c_{h(n)} + \operatorname{sgn}(\gamma_{h(n)}) \bigl( |\gamma_{h(n)} | -
\beta\bigl(\varepsilon,h(n)\bigr) \bigr)_{+} ;
\]
\end{itemize}
until $\textit{best}\_\textit{improvement}<\kappa$.

\bigskip

Return the estimation:
\[
[\tilde{Y}_{N+1},\ldots,\tilde{Y}_{(k+1)N} ] = \bigl[\hat{\theta
}(X_{N+1}),\ldots,\hat{\theta}(X_{(k+1)N}) \bigr],
\]
where:
\[
\hat{\theta}(\cdot) = \sum_{h=1}^{m} c_{h} \theta_{h} (\cdot) .
\]
\\
\hline
\end{tabular}
\caption{Detailed version of the feature selection algorithm in the
transductive case.}
\end{figure}

\begin{TEO} We have:
\[
P_{2N} \bigl[\forall n\in\{1,\ldots,n_{0}\}, r_{2} \bigl(\theta^{ (n )} \bigr)\leq
r_{2} \bigl(\theta^{ (n-1 )} \bigr)
-d_{2}^{2} \bigl(\theta^{(n)},\theta^{(n-1)} \bigr) \bigr] \geq
1-\varepsilon.
\]
\end{TEO}

The proof of this theorem is exactly the same as the proof of
Theorem \ref{propalgo}.

\begin{exm}[(Estimation of wavelet coefficients)]
Let us consider the case where $\varTheta_{0}$ does not depend on the observations.
We can, for example, choose a basis of $\varTheta$, or a basis of a
subspace of $\varTheta$.
We obtain an estimator of the form:
\[
\theta(x)=\sum_{h=1}^{m} \alpha^{h} \theta_{h}(x) .
\]
In the case when $(\theta_{k})_{k}$ is a wavelet basis, then we
obtain here again a procedure for thresholding wavelets
coefficients.
\end{exm}

\begin{exm}[(SVM and multiscale SVM)]
Let us choose $\varTheta$ as the set of all functions
$\mathcal{X}\rightarrow\mathdsR $, a family of kernels
$K_{1},\ldots,K_{m'(N)}$ for a $m'(N)\geq1$ and:
\[
\varTheta_{0} = \bigl\{ K_{h}(X_{i},\cdot) , h\in\bigl\{1,\ldots,m'(N)\bigr\}, i\in\bigl\{1,\ldots
,(k+1)N\bigr\} \bigr\}.
\]
In this case we have $m=(k+1)Nm'(N)$.
We obtain an estimator of the form:
\[
\theta(x)=\sum_{h=1}^{m'(N)} \sum_{j=1}^{2N} \alpha^{j,h}
K_{h}(X_{j},x) .
\]
Let us put:
\[
I_{h}= \bigl\{j\in\{1,\ldots,2N\},\alpha^{j,h}\neq0 \bigr\} .
\]
We have:
\[
\theta(x)=\sum_{h=1}^{m'(N)} \sum_{j\in I_{h}} \alpha^{j,h}
K_{h}(X_{i},x) ,
\]
that is a Support Vector Machine with different kernel estimate;
like in Example \ref{MSVM}, the kernels $K_{h}$ can be the same
kernel taken at different scales.
\end{exm}

\begin{exm}[(Kernel PCA Kernel Projection Machine)]
Take the same $\varTheta$ and consider the kernel:
\[
K\bigl(x,x'\bigr) = \bigl\<\varPsi(x),\varPsi\bigl(x'\bigr) \bigr\> .
\]
Let us consider a principal component analysis (PCA) of the family:
\[
\bigl\{K(X_{i},\cdot),\ldots,K(X_{(k+1)N},\cdot) \bigr\}
\]
by performing a diagonalization of the matrix:
\[
\bigl(K (X_{i},X_{j} ) \bigr)_{1\leq i,j \leq(k+1)N} .
\]
This method is known as Kernel PCA, see for example~\cite{PCA}.
We obtain eigenvalues:
\[
\lambda^{1} \geq\cdots\geq\lambda^{(k+1)N}
\]
and associated eigenvectors $ e^{1},\ldots, e^{(k+1)N} $, associated
to elements of $\varTheta$:
\[
k_{1}(\cdot) = \sum_{i=1}^{(k+1)N} e^{1}_{i}K(X_{i},\cdot) ,\ldots,
k_{(k+1)N}(\cdot) = \sum_{i=1}^{(k+1)N} e^{(k+1)N}_{i}K(X_{i},\cdot)
\]
that are exchangeable functions of the observations.
Using the family:
\[
\varTheta_{0} = \{k_{1},\ldots,k_{(k+1)N} \},
\]
we obtain an algorithm that selects which eigenvectors are going to be used
in the regression estimation. This is very close to the Kernel
Projection Machine (KPM)
described by Blanchard, Massart, Vert and Zwald \cite{KPM} in the
context of
classification.
\end{exm}

\section{Rates of convergence in Sobolev and Besov spaces}

\label{rate}

We conclude this paper by coming back to the inductive case. We use
Theorem \ref{propalgo} as an oracle inequality to show that the
obtained estimator is adaptative, which means that if we assume that
the true regression function $f$ has an unknown regularity $\beta$,
then the estimator is able to reach the optimal speed of convergence
$N^{{-2\beta}/(2\beta+1)}$ up to a $\log N$ factor.

\subsection{Presentation of the context}

Here we assume that $\mathcal{X}$ is a compact interval of
$\mathdsR $, that $\varTheta=\mathdsL_{2}(P_{(X)})$ and that $P$ is
such that $Y=f(X)+\eta$ with $\eta$ independent of $X$, $P\eta=0$
and $P(\eta^{2})\leq\sigma^{2} < +\infty$.

We assume that $(\theta_{k})_{k\in\mathdsN ^{*}}$ is an orthonormal
basis of $\varTheta$. We still have to choose $m\in\mathdsN $ and we
will take $\varTheta_{0}=\{\theta_{1},\ldots,\theta_{m}\}$.

Remark that the orthogonality means here that
$P[\theta_{k}(X)^{2}]=1$ for any $k\in\mathdsN ^{*}$, and that:
\[
P\bigl[\theta_{k}(X)\theta_{k'}(X)\bigr]=0
\]
for any $k'\neq k$.

\subsection{Rate of convergence of the estimator: the Sobolev space case}

Now, let us put:
\[
\overline{\theta}^{m}=\mathop{\arg\min}_{\theta\in \operatorname{Span}(\varTheta_{0})} R(\theta)
\]
(that depends effectively on $m$ by $\varTheta_{0}=\{\theta_{1},\ldots
,\theta_{m}\}$),
and let us assume that $f$
satisfies the two following conditions: it is regular, namely there is
an unknown $\beta\geq1$ and a $C\geq0$ such that:
\[
\big\|\overline{\theta}^{m}-f \big\|_{P}^{2} \leq C m^{-2\beta},
\]
and that we have a constant $B<\infty$ such that:
\[
\sup_{x\in\mathcal{X}} f(x) \leq B
\]
with $B$ known to the statistician.
It follows that:
\[
\|f \|^{2}_{P} \leq B^{2} .
\]
If follows that every set, for $k\in\{1,\ldots,m\}$:
\[
\mathcal{F}_{k} = \Biggl\{\sum_{j=1}^{\infty}\alpha_{j}\theta_{j}\dvt  \alpha
_{k}^{2}\leq B^{2} \Biggr\} \cap\varTheta
\]
is a convex set that contains $f$ and such that the orthogonal
projection: $\varPi_{P}^{\mathcal{F},m} =
\varPi_{P}^{\mathcal{F}_{m}}\cdots\varPi_{P}^{\mathcal{F}_{1}}$ (where
$\varPi_{P}^{\mathcal{F}_{k}}$ denotes the orthogonal projection on
$\mathcal{F}_{k}$) can only improve an estimator:
\[
\forall\theta,\quad \big\|\varPi_{P}^{\mathcal{F},m}\theta - f \big\|_{P}^{2}\leq \|
\theta - f \|_{P}^{2}.
\]
Actually, note that this projection just consists in thresholding very
large coefficients to a limited value.
This modification is necessary in what follows, but this is just a
technical remark: most of the time, our estimator
won't be modified by $\varPi_{P}^{\mathcal{F},m}$ for any $m$.

Remember also that in this context, the estimator given in
Definition \ref{hatf} is just:
\[
\hat{\theta} = \varPi_{P}^{m,\varepsilon} \cdots\varPi
_{P}^{1,\varepsilon} 0 .
\]

\begin{TEO} \label{speed}
Let us assume that $\varTheta=\mathdsL_{2}(P_{(X)})$, $\mathcal{X}=
[0,1]$ and $(\theta_{k})_{k\in\mathdsN ^{*}}$
is an orthonormal basis of $\varTheta$.
Let us assume that we are in the idealized regression model:
\[
Y=f(X)+\eta,
\]
where $P\eta=0$, $P(\eta^{2})\leq\sigma^{2}<\infty$ and $\eta$ and
$X$ are independent, and $\sigma$ is known. Let us assume that
$f\in\varTheta$ is such that there is an unknown $\beta\geq1$ and an
unknown $C\geq0$ such that:
\[
\|\overline{\theta}_{m}-f \|_{P}^{2} \leq C m^{-2\beta},
\]
and that we have a constant $B<\infty$ such that:
\[
\sup_{x\in\mathcal{X}} f(x) \leq B
\]
with $B$ known to the statistician. Then our estimator
$\hat{\theta}$ (given in Definition \textup{\ref{hatf}} with $n_{0}=m$ here,
build using the bound $\beta(\varepsilon,k)$ given in Theorem
\textup{\ref{lastTH}}), with $\varepsilon=N^{-2}$ and $m=N$, is such that,
for any $N\geq2$,
\[
P^{\otimes N} \bigl[ \big\|\varPi_{P}^{\mathcal{F},N}\hat{\theta}-f \big\|_{P}^{2} \bigr]
\leq C'(C,B,\sigma) \biggl(\frac{\log N}{N} \biggr)^{{2\beta}/{(2\beta+1)}} .
\]
\end{TEO}

Here again, the proof is given at the end of the paper
(Section~\ref{subproofrate}). Let us just remark
that, in the case where $\mathcal{X}=[0,1]$, $P$ is the Lebesgue
measure, and $(\theta_{k})_{k\in\mathdsN ^{*}}$ is the
trigonometric basis, the condition:
\[
\big\|\overline{\theta}^{m}-f \big\|_{P}^{2} \leq C m^{-2\beta}
\]
is satisfied for $C=C(\beta,L)$ as soon as $f\in W(\beta,L)$ where
$W(\beta,L)$ is the Sobolev class:
\[
\biggl\{f\in\mathcal{L}^{2}\dvt  f^{(\beta-1)} \mbox{ is absolutely continuous
and }
\int_{0}^{1}f^{(\beta)}(x)^{2}\lambda(\mathrm{d}x) \leq L^{2}
\biggr\}.
\]
The minimax rate of convergence in $W(\beta,L)$ is $N^{-{2\beta
}/{(2\beta+1)}}$,
so we can see that our estimator reaches the best rate of convergence
up to a $\log N$ factor
with an unknown $\beta$.

\subsection{Rate of convergence in Besov spaces}

We here extend the previous result to the case of a Besov space
$B_{s,p,q}$ in the case of a wavelet basis
(see \cite{Ondel2} or~\cite{Ondel}).

\begin{TEO} \label{rate2}
Let us assume that $\mathcal{X}=[-A,A]$, that $P_{(X)}$ is uniform on
$\mathcal{X}$ and that
$(\psi_{j,k})_{j=0,\ldots,+\infty, k\in\{1,\ldots,2^{j}\}}$ is a
wavelet basis, together with
a function $\phi$, satisfying the conditions given in \cite{Ondel2},
with $\phi$ and $\psi_{0,1}$ supported by $[-A,A]$.
Let us assume that $f\in B_{s,p,q}$ with $s>\frac{1}{p}$, $1\leq
p,q\leq\infty$, with:
\begin{eqnarray*}
B_{s,p,q} &=&  \Biggl\{g\dvt [-A,A]\rightarrow\mathdsR ,
g(\cdot)=\alpha\phi(\cdot) + \sum_{j=0}^{\infty}\sum_{k=1}^{2^{j}}\beta
_{j,k}\psi_{j,k}(\cdot),
\\
&& \sum_{j=0}^{\infty}2^{jq (s-{1}/{2}-{1}/{p} )} \Biggl[\sum
_{k=1}^{2^{j}} |\beta_{j,k} |^{p} \Biggr]^{{q}/{p}}
= \|g\|_{s,p,q}^{q} < + \infty
\Biggr\}
\end{eqnarray*}
(with obvious changes for $p=+\infty$ or $q=+\infty$) with unknown
constants $s$, $p$ and $q$ and that for any $x$, $|f(x)|\leq B$ for
a known constant $B$. Let us choose:
\[
\{\theta_{1},\ldots,\theta_{m}\} = \{\phi\} \cup\bigl\{\psi_{j,k},j=1,\ldots
,2^{\lfloor {\log N}/{\log2}\rfloor},k=1,\ldots,2^{j}\bigr\}
\]
(so $\frac{N}{2} \leq m\leq N$) and $\varepsilon=N^{-2}$ in the
definition of $\hat{\theta}$. Then we have:
\[
P^{\otimes N} \bigl[ \big\|\varPi_{P}^{\mathcal{F},N}\hat{\theta}-f \big\|_{P}^{2} \bigr] =
\EuScript{O} \biggl( \biggl(\frac{\log N}{N} \biggr)^{{2s}/{(2s+1)}}
(\log N )^{ (1-{2}/{((1+2s)q)} )_{+}}
\biggr).
\]
\end{TEO}

Let us remark that we obtain nearly the same rate of convergence
than in \cite{Ondel2}, namely the minimax rate of convergence up to
a $\log N$ factor.

For the proof, see Section~\ref{subproofrate}.

\section{Proofs}\label{proofs}

The order of the proofs is exactly the order of apparition of the
results in the paper, except for the first theorem (Theorem
\ref{lastTH}): its proof using lemmas proved in the transductive
setting, it is given after the proof of the transductive theorems.

\subsection[Proof of Theorems 2.4--2.6]{Proof of Theorems \textup{{\protect\ref{TH1}--\protect\ref{THSVM}}}}
\label{subproof1}

First, we prove a lemma that is the basis of proofs of Theorems
\ref{TH1}--\ref{THSVM}.

\begin{lemma} \label{inequality}
We have, for any $\theta\in\varTheta$, $\gamma> 0$ and $\eta\geq0$:
\[
P\exp(\gamma W_{\theta} -\eta)  =
\exp\biggl\{\frac{\gamma^{2}}{2}V (W_{\theta} )
+\frac{\gamma^{3}}{2}\int_{0}^{1}(1-\beta)^{2}M^{3}_{\gamma\beta
W_{\theta}} (W_{\theta} ) \,\mathrm{d}\beta-\eta\biggr\} ,
\]
and
\[
P\exp(-\gamma W_{\theta} -\eta)  =
\exp\biggl\{\frac{\gamma^{2}}{2}V (W_{\theta} )
-\frac{\gamma^{3}}{2}\int_{0}^{1}(1-\beta)^{2}M^{3}_{\gamma\beta
W_{\theta}} (W_{\theta} ) \,\mathrm{d}\beta-\eta\biggr\}.
\]
\end{lemma}

\begin{pf}
For the first equality, we write:
\begin{eqnarray*}
\log P\exp(\gamma W_{\theta} -\eta)
&=&  \log P \exp (\gamma W_{\theta} ) - \eta
\\
&=&  \int_{0}^{\gamma} P_{\beta
W_{\theta}} (W_{\theta} )\,\mathrm{d}\beta-\eta= \int_{0}^{\gamma}
(\gamma-\beta) V_{\beta W_{\theta}} (W_{\theta} ) \,\mathrm{d}\beta-
\eta
\\
&=&  \frac{\gamma^{2}}{2}V (W_{\theta} ) +
\int_{0}^{\gamma}\frac{(\gamma-\beta)^{2}}{2} M^{3}_{\beta
W_{\theta}}(W_{\theta} ) \,\mathrm{d}\beta- \eta
\\
&=&  \frac{\gamma^{2}}{2}V (W_{\theta} )+\frac{\gamma^{3}}{2}
\int_{0}^{1}(1-\beta)^{2} M^{3}_{\gamma\beta W_{\theta}}
(W_{\theta} ) \,\mathrm{d}\beta- \eta.
\end{eqnarray*}
For the reverse equality, the proof is exactly the same, replacing
$\gamma$ by $-\gamma$.
\end{pf}

We can now give the proof of both theorems.

\begin{pf*}{Proof of Theorem \ref{TH1}}
Let us choose $k\in\{1,\ldots,m\}$, for any $\lambda_{k}> 0$ and
$\eta_{k} \geq0$ we have:
\begin{eqnarray*}
\lefteqn{P^{\otimes N} \exp \Biggl\{\frac{\lambda_{k}}{N}\sum_{i=1}^{N}
\bigl[Y_{i}\theta_{k}(X_{i})-P \bigl(Y\theta_{k}(X) \bigr) \bigr]
-\eta_{k} \Biggr\}}
\\
\lefteqn{\quad =  \biggl\{ P \exp \biggl[\frac{\lambda_{k}}{N} W_{\theta_{k}} - \frac{\eta
_{k}}{N} \biggr] \biggr\}^{N}}
\\
\lefteqn{\quad =  \exp\biggl[ \frac{\lambda_{k}^{2}}{2N}V (W_{\theta_{k}} )
+\frac{\lambda_{k}^{3}}{2N^{2}}\int_{0}^{1}(1-\beta)^{2}M^{3}_{
({\beta\lambda_{k}}/{N}) W_{\theta_{k}}} (W_{\theta_{k}} ) \,\mathrm{d}\beta
- \eta_{k} \biggr]}
\end{eqnarray*}
by the first equality of Lemma \ref{inequality}. By the same way,
using the reverse inequality we obtain:
\begin{eqnarray*}
\lefteqn{P^{\otimes N} \exp \Biggl\{\frac{\lambda_{k}}{N}\sum_{i=1}^{N}
\bigl[P \bigl(Y\theta_{k}(X) \bigr)-Y_{i}\theta_{k}(X_{i}) \bigr]
-\eta_{k} \Biggr\}}
\\
\lefteqn{\quad = \exp\biggl[ \frac{\lambda_{k}^{2}}{2N}V (W_{\theta_{k}} )
-\frac{\lambda_{k}^{3}}{2N^{2}}\int_{0}^{1}(1-\beta)^{2}M^{3}_{
(\beta\lambda_{k}/N) W_{\theta_{k}}} (W_{\theta_{k}} )\,\mathrm{d}\beta
- \eta_{k} \biggr].}
\end{eqnarray*}
So we obtain, for any $k\in\{1,\ldots,m\}$, for any $\lambda_{k}> 0$
and $\eta_{k} \geq0$:
\begin{eqnarray*}
\lefteqn{P^{\otimes N} \exp \Biggl\{\lambda_{k} \Bigg|
\frac{1}{N}\sum_{i=1}^{N}Y_{i}\theta_{k}(X_{i})-P \bigl(Y\theta_{k}(X) \bigr) \Bigg|
-\eta_{k} \Biggr\}}
\\
\lefteqn{\quad \leq 2 \exp\biggl[
\frac{\lambda_{k}^{2}}{2N}V (W_{\theta_{k}} )-\eta_{k} \biggr]
\cosh\biggl[
\frac{\lambda_{k}^{3}}{2N^{2}}\int_{0}^{1}(1-\beta)^{2}M^{3}_{
(\beta\lambda_{k}/N) W_{\theta_{k}}} (W_{\theta_{k}} )\,\mathrm{d}\beta
\biggr]}
\\
\lefteqn{\quad \leq 2 \exp\biggl[
\frac{\lambda_{k}^{2}}{2N}V (W_{\theta_{k}} )-\eta_{k}
+ \frac{\lambda_{k}^{6}}{8N^{4}} \biggl(\int_{0}^{1}(1-\beta)^{2}M^{3}_{
(\beta\lambda_{k}/N) W_{\theta_{k}}} (W_{\theta_{k}} )\,\mathrm{d}\beta\biggr)^{2}
\biggr] ,}
\end{eqnarray*}
since, for any $x\in\mathdsR $, we have:
\[
\cosh(x) \leq\exp\biggl(\frac{x^{2}}{2} \biggr).
\]
Now, let us choose $\varepsilon>0$ and put:
\[
\eta_{k}= \frac{\lambda_{k}^{2}}{2N}V (W_{\theta_{k}} )
+ \frac{\lambda_{k}^{6}}{8N^{4}} \biggl(\int_{0}^{1}(1-\beta)^{2}M^{3}_{
({\beta\lambda_{k}}/{N}) W_{\theta_{k}}} (W_{\theta_{k}} )\,\mathrm{d}\beta\biggr)^{2}
-\log\frac{\varepsilon}{2m}.
\]
We obtain:
\begin{eqnarray*}
\lefteqn{P^{\otimes N} \sum_{k=1}^{m} \exp \Biggl\{\lambda_{k} \Bigg|
\frac{1}{N}
\sum_{i=1}^{N}Y_{i}\theta_{k}(X_{i})-P \bigl(Y\theta_{k}(X) \bigr)
\Bigg|}
\\
\lefteqn{\quad {}
-
\frac{\lambda_{k}^{2}}{2N}V (W_{\theta_{k}} )
+ \frac{\lambda_{k}^{6}}{8N^{4}} \biggl(\int_{0}^{1}(1-\beta)^{2}M^{3}_{
({\beta\lambda_{k}}/{N}) W_{\theta_{k}}} (W_{\theta_{k}} )\,\mathrm{d}\beta\biggr)^{2}
+\log\frac{\varepsilon}{2m} \Biggr\}
\leq\varepsilon}
\end{eqnarray*}
and so:
\begin{eqnarray*}
\lefteqn{P^{\otimes N} \Biggl[ \forall k \in\{1,\ldots,m\},   \Bigg|
\frac{1}{N}
\sum_{i=1}^{N}Y_{i}\theta_{k}(X_{i})-P \bigl(Y\theta_{k}(X) \bigr)
\Bigg|}
\\
\lefteqn{\quad \leq\frac{\lambda_{k}}{2N}V (W_{\theta_{k}} )
+ \frac{\lambda_{k}^{5}}{8N^{4}} \biggl(\int_{0}^{1}(1-\beta)^{2}M^{3}_{
({\beta\lambda_{k}}/{N}) W_{\theta_{k}}} (W_{\theta_{k}} )\,\mathrm{d}\beta\biggr)^{2}
+ \frac{\log ({2m}/{\varepsilon})}{\lambda_{k}} \Biggr]
\geq1 - \varepsilon.}
\end{eqnarray*}
Now, we put:
\[
\lambda_{k} = \sqrt{\frac{2N\log ({2m}/{\varepsilon}) }{V (W_{\theta
_{k}} )}} .
\]
We obtain, with $P^{\otimes N}$-probability at least
$1-\varepsilon$, for any $k\in\{1,\ldots,m\}$:
\begin{eqnarray*}
\Bigg| \frac{1}{N}
\sum_{i=1}^{N}Y_{i}\theta_{k}(X_{i})-P \bigl(Y\theta_{k}(X) \bigr) \Bigg|
&\leq& \sqrt{\frac{2V (W_{\theta_{k}} )\log
({2m}/{\varepsilon})}{N}}
\\
&&{} + \frac{\log^{{5}/{2}} ({2m}/{\varepsilon})}{N V (W_{\theta_{k}} )^{3}}
\biggl(\int_{0}^{1}(1-\beta)^{2}M^{3}_{({\beta\lambda_{k}}/{N}) W_{\theta
_{k}}} (W_{\theta_{k}} )\,\mathrm{d}\beta\biggr)^{2}.
\end{eqnarray*}
For short, we take the notation of the theorem:
\[
I_{\theta_{k}} (\gamma ) = \int_{0}^{1}(1-\beta)^{2}M^{3}_{\beta
\gamma W_{\theta_{k}}} (W_{\theta_{k}} ) .
\]
Now, dividing both sides by:
\[
P \bigl[\theta_{k}(X)^{2} \bigr]
\]
we obtain:
\[
| \hat{\alpha}_{k}\mathcal{C}_{k} -
\overline{\alpha}_{k} | \leq
\frac{1}{P [\theta_{k}(X)^{2} ]} \biggl[
\sqrt{\frac{2V (W_{\theta_{k}} )\log
({2m}/{\varepsilon})}{N}}
+ \frac{I_{\theta_{k}}^{2} ({\lambda_{k}}/{N} ) \log^{{5}/{2}}
({2m}/{\varepsilon})}{N V (W_{\theta_{k}} )^{3}} \biggr].
\]
In order to conclude, just remark that:
\[
R(\hat{\alpha}_{k}\mathcal{C}_{k}\theta_{k})-R(\overline{\alpha
}_{k}\theta_{k})
= | \hat{\alpha}_{k}\mathcal{C}_{k} -
\overline{\alpha}_{k} |^{2} P \bigl[\theta_{k}(X)^{2} \bigr] .
\]\upqed
\end{pf*}

\begin{pf*}{Proof of Theorem \ref{TH1bis}}
Remark that, for any $\theta\in\varTheta$:
\[
V (W_{\theta} ) = P \bigl(W_{\theta}^{2} \bigr) - P (W_{\theta} )^{2} ,
\]
we will deal with each term separately. For the first term, let us
remark that we obtain the following result that is obtained exactly
as Lemma \ref{inequality}. For any $\theta\in\varTheta$:
\begin{eqnarray*}
P\exp\bigl\{\gamma \bigl[P \bigl(W_{\theta}^{2} \bigr)
-W_{\theta}^{2} \bigr] -\eta\bigr\}
= \exp\biggl\{\frac{\gamma^{2}}{2}V \bigl(W_{\theta}^{2} \bigr)
+\frac{\gamma^{3}}{2}\int_{0}^{1}(1-\beta)^{2}M^{3}_{\gamma\beta
W_{\theta}^{2}} \bigl(W_{\theta}^{2} \bigr)\,\mathrm{d}\beta-\eta\biggr\}.
\end{eqnarray*}
Let us apply this result to every $\theta_{k}$ for
$k\in\{1,\ldots,m\}$:
\begin{eqnarray*}
P^{\otimes N}\exp \Biggl\{\lambda_{k}
\Biggl[P \bigl(W_{\theta_{k}}^{2} \bigr)
-\frac{1}{N}\sum_{i=1}^{N}Y_{i}^{2}\theta_{k}(X_{i})^{2} \Biggr]
-\eta_{k} \Biggr\}
= \exp\biggl\{\frac{\lambda_{k}^{2}}{2N}V \bigl(W_{\theta_{k}}^{2} \bigr)
+\frac{\lambda_{k}^{3}}{2N} J_{k} \biggl(\frac{\lambda_{k}}{N} \biggr) -\eta_{k} \biggr\},
\end{eqnarray*}
where:
\[
J_{\theta} (\gamma) = \int_{0}^{1}(1-\beta)^{2}M^{3}_{\gamma\beta
W_{\theta_{k}}^{2}} \bigl(W_{\theta}^{2} \bigr)\,\mathrm{d}\beta.
\]
Taking
\[
\eta_{k} = \frac{\lambda_{k}^{2}}{2N}V \bigl(W_{\theta_{k}}^{2} \bigr)
+\frac{\lambda_{k}^{3}}{2N^{2}} J_{\theta_{k}} \biggl(\frac{\lambda_{k}}{N}
\biggr) + \log\frac{2m}{\varepsilon}
\]
and
\[
\lambda_{k} = \sqrt{\frac{2N\log ({2m}/{\varepsilon})}{V (W_{\theta
_{k}}^{2} )}}
\]
we obtain that the following inequality is satisfied with
$P^{\otimes N}$-probability at least $1-\frac{\varepsilon}{2}$, for
any $k$:
\begin{eqnarray} \label{inequ1}
\nonumber P \bigl(W_{\theta_{k}}^{2} \bigr) &\leq&
\frac{1}{N}\sum_{i=1}^{N}Y_{i}^{2}\theta_{k}(X_{i})^{2} +
\sqrt{\frac{2V (W_{\theta_{k}}^{2} )\log ({2m}/{\varepsilon})}{N}}
+ \frac{\log ({2m}/{\varepsilon})}{N
V (W_{\theta_{k}}^{2} ) }
J_{\theta_{k}} \biggl(\sqrt{\frac{2\log ({2m}/{\varepsilon})}{NV (W_{\theta
_{k}}^{2} )}} \biggr)
\\
&=&  \frac{1}{N}\sum_{i=1}^{N}Y_{i}^{2}\theta_{k}(X_{i})^{2} +
\mathcal{A}_{k}
\end{eqnarray}
for short. Now, we try to upper bound the second term,
$-P (W_{\theta} )^{2}$. Remark that, for any $\theta$:
\begin{eqnarray*}
\Biggl(\frac{1}{N}\sum_{i=1}^{N} Y_{i} \theta(X_{i}) \Biggr)^{2} -
P (W_{\theta} )^{2}
&=&  \Biggl(\frac{1}{N}\sum_{i=1}^{N} Y_{i} \theta(X_{i}) - P (W_{\theta} ) \Biggr)
\Biggl(\frac{1}{N}\sum_{i=1}^{N} Y_{i} \theta(X_{i}) + P (W_{\theta} ) \Biggr)
\\
&\leq& \Bigg|\frac{1}{N}\sum_{i=1}^{N} Y_{i} \theta(X_{i}) - P (W_{\theta} ) \Bigg|
\\
&&{}\times \Biggl\{ 2 \Bigg|\frac{1}{N}\sum_{i=1}^{N} Y_{i} \theta(X_{i}) \Bigg|
+ \Bigg|\frac{1}{N}\sum_{i=1}^{N} Y_{i} \theta(X_{i}) - P (W_{\theta} ) \Bigg| \Biggr\}.
\end{eqnarray*}
Remember that in the proof of Theorem \ref{TH1} we got the upper
bound, with probability at least $1-\frac{\varepsilon}{2}$, for any
$k$:
\begin{eqnarray*}
\Bigg| \frac{1}{N}
\sum_{i=1}^{N}Y_{i}\theta_{k}(X_{i})-P \bigl(Y\theta_{k}(X) \bigr) \Bigg|
\leq\sqrt{\frac{2V (W_{\theta_{k}} )\log
({4m}/{\varepsilon})}{N}}
+ \frac{\log^{{5}/{2}} ({4m}/{\varepsilon})}{N V (W_{\theta_{k}} )^{3}}
I_{\theta_{k}} \biggl(\sqrt{\frac{2\log ({4m}/{\varepsilon})}{N V (W_{\theta
_{k}} )}} \biggr)^{2}
,
\end{eqnarray*}
that gives:
\begin{eqnarray} \label{inequ2}
\nonumber - P (W_{\theta_{k}} )^{2} &\leq& -
\Biggl(\frac{1}{N}\sum_{i=1}^{N} Y_{i} \theta_{k}(X_{i}) \Biggr)^{2}
+ \biggl\{ \sqrt{\frac{2V (W_{\theta_{k}} )\log
({4m}/{\varepsilon})}{N}}
+ \frac{\log^{{5}/{2}} ({4m}/{\varepsilon})}{N V (W_{\theta_{k}} )^{3}}
I_{\theta_{k}} \biggl(\sqrt{\frac{2\log ({4m}/{\varepsilon})}{N V (W_{\theta
_{k}} )}} \biggr)^{2} \biggr\}
\\
\nonumber &&{}\times \Biggl\{ 2 \Bigg|\frac{1}{N}\sum_{i=1}^{N} Y_{i} \theta_{k}(X_{i})
\Bigg| +
\sqrt{\frac{2V (W_{\theta_{k}} )\log
({4m}/{\varepsilon})}{N}}
+ \frac{\log^{{5}/{2}} ({4m}/{\varepsilon})}{N V (W_{\theta_{k}} )^{3}}
I_{\theta_{k}} \biggl(\sqrt{\frac{2\log ({4m}/{\varepsilon})}{N V (W_{\theta
_{k}} )}} \biggr)^{2}
\Biggr\}
\\
&=&  - \Biggl(\frac{1}{N}\sum_{i=1}^{N} Y_{i} \theta(X_{i}) \Biggr)^{2} + \mathcal{B}_{k}
\end{eqnarray}
for short. Let us combine inequalities (\ref{inequ1})~and~(\ref{inequ2}).
We obtain that, with probability at least
$1-\varepsilon$, for every $k$ we have:
\[
V (W_{\theta_{k}} ) =  P \bigl(W_{\theta_{k}}^{2} \bigr) -
P (W_{\theta_{k}} )^{2}
\leq\frac{1}{N}\sum_{i=1}^{N} Y_{i}^{2} \theta_{k}(X_{i})^{2}-
\Biggl(\frac{1}{N}\sum_{i=1}^{N} Y_{i} \theta_{k}(X_{i}) \Biggr)^{2}
+\mathcal{A}_{k} +\mathcal{B}_{k}
=  \hat{V}_{k} +\mathcal{A}_{k} +\mathcal{B}_{k}.
\]\upqed
\end{pf*}

\begin{pf*}{Proof of Theorem \ref{THSVM}}
This proof is a variant of the proof of Theorem \ref{TH1}, the
method it uses is due to Seeger \cite{Seeger}. Let us define, for
any $i\in\{1,\ldots,N\}$:
\[
P_{i}(\cdot) = P^{\otimes N} (\cdot |Z_{i}) .
\]
Let us choose $(i,k)\in\{1,\ldots,N\}\times\{1,\ldots,m'(N)\}$, for any
$\lambda_{i,k}=\lambda_{i,k}(Z_{i})> 0$ and
$\eta_{i,k}=\eta_{i,k}(Z_{i}) \geq0$ we have:
\begin{eqnarray*}
\lefteqn{P_{i} \exp \biggl\{\frac{\lambda_{i,k}}{N-1}\sum_{j\neq i}
\bigl[Y_{j}\theta_{i,k}(X_{j})-P \bigl(Y\theta_{i,k}(X) \bigr)
\bigr] -\eta_{i,k} \biggr\}}
\\
\lefteqn{\quad \leq\exp\Biggl[ \frac{\lambda_{i,k}}{2(N-1)} V (W_{\theta_{i,k}} )
+ \frac{\lambda_{i,k}^{3}}{2(N-1)^{2}} \int_{0}^{1} (1-\beta
)^{2}M^{3}_{(\beta\lambda_{i,k}/{N-1})
W_{\theta_{i,k}}} (W_{\theta_{i,k}} )\,\mathrm{d}\beta
-\eta_{i,k} \Biggr]}
\end{eqnarray*}
by the first equality of Lemma \ref{inequality}. In the same way, we
obtain the reverse inequality and, combining both results, for any
$(i,k)\in\{1,\ldots,N\}\times\{1,\ldots,m'(N)\}$, for any $\lambda_{i,k}>
0$ and $\eta_{i,k} \geq0$:
\begin{eqnarray*}
\lefteqn{P_{i} \exp \biggl\{\lambda_{i,k} \bigg| \frac{1}{N-1}\sum_{j\neq
i}Y_{j}\theta_{i,k}(X_{j})
-P \bigl(Y\theta_{i,k}(X) \bigr) \bigg|-\eta_{i,k} \biggr\}}
\\
\lefteqn{\quad \leq2 \exp\biggl[ \frac{\lambda_{i,k}^{2}}{2(N-1)}V (W_{\theta_{i,k}}
)-\eta_{i,k} \biggr] \cosh\biggl[
\frac{\lambda_{i,k}^{3}}{2(N-1)^{2}}I_{i,k} \biggr]}
\\
\lefteqn{\quad \leq2 \exp\biggl[
\frac{\lambda_{i,k}^{2}}{2(N-1)}V (W_{\theta_{i,k}} )-\eta_{i,k}
+ \frac{\lambda_{i,k}^{6}}{8(N-1)^{4}}I_{i,k}^{2} \biggr],}
\end{eqnarray*}
where:
\[
I_{ i,k} =
\int_{0}^{1}(1-\beta)^{2}M^{3}_{{(\beta\lambda_{i,k}/{N})}
W_{\theta_{i,k}}} (W_{\theta_{i,k}} )\,\mathrm{d}\beta
\]
for short.
Now, let us choose $\varepsilon>0$ and put:
\[
\eta_{i,k}= \frac{\lambda_{i,k}^{2}}{2(N-1)}V (W_{\theta_{i,k}} )
+ \frac{\lambda_{i,k}^{6}}{8(N-1)^{4}} I_{i,k}^{2}
-\log\frac{\varepsilon}{2Nm'(N)}.
\]
We obtain:
\begin{eqnarray*}
\lefteqn{P^{\otimes N} \sum_{i=1}^{N} \sum_{k'=1}^{m'(N)}
\exp \biggl\{\lambda_{i,k} \bigg| \frac{1}{N-1} \sum_{j\neq i} Y_{j}\theta
_{i,k}(X_{j})-P \bigl(Y\theta_{i,k}(X) \bigr) \bigg|}
\\
\lefteqn{\qquad{}
- \frac{\lambda_{i,k}^{2}}{2(N-1)}V (W_{\theta_{i,k}} )
- \frac{\lambda_{i,k}^{6}}{8(N-1)^{4}}I_{i,k}^{2} + \log\frac
{\varepsilon}{2Nm'(N)}
\biggr\}}
\\
\lefteqn{\quad = P^{\otimes N} \sum_{i=1}^{N} \sum_{k'=1}^{m'(N)} P_{i}
\exp \biggl\{\lambda_{i,k} \bigg| \frac{1}{N-1} \sum_{j\neq i} Y_{j}\theta
_{i,k}(X_{j})-P \bigl(Y\theta_{i,k}(X) \bigr) \bigg|}
\\
\lefteqn{\qquad{}
- \frac{\lambda_{i,k}^{2}}{2(N-1)}V (W_{\theta_{i,k}} )
- \frac{\lambda_{i,k}^{6}}{8(N-1)^{4}}I_{i,k}^{2} + \log\frac
{\varepsilon}{2Nm'(N)}
\biggr\}
\leq\varepsilon.}
\end{eqnarray*}
Now, we put:
\[
\lambda_{i,k} = \sqrt{\frac{2N\log ({2Nm'(N)}/{\varepsilon}) }{V
(W_{\theta_{i,k}} )}} ,
\]
and achieve the proof exactly as for Theorem \ref{TH1}.
\end{pf*}

\subsection[Proof of Theorems 3.1 and 3.2]{Proof of Theorems \textup{{\protect\ref{thtrans}}} and \textup{{\protect\ref{cortrans}}}}

\label{subproof2}

Here again, the first thing to do is to prove a general deviation
inequality. This one is a variant of the one given by Catoni
\cite{Classif}. We go back to the notations of Theorem \ref{thtrans}
and \ref{cortrans}, with test sample of size $N$.

\begin{dfn}
Let $\mathcal{G}$ denote the set of all functions:
\begin{eqnarray*}
\lefteqn{g \dvtx  (\mathcal{X}\times\mathdsR  )^{2N} \times\mathdsR ^{2}
\rightarrow\mathdsR ,}
\\
\lefteqn{\bigl(Z_{1},\ldots,Z_{2N},u,u' \bigr)  \mapsto g \bigl(Z_{1},\ldots,Z_{2N},u,u' \bigr) =
g\bigl(u,u'\bigr)}
\end{eqnarray*}
for the sake of simplicity, such that $g$ is exchangeable with
respect to its $2N$ first arguments.
\end{dfn}

\begin{lemma} \label{learninglemma}
For any exchangeable probability distribution $\mathcal{P}$ on
$ (Z_{1},\ldots,Z_{2N} )$, for any measurable function
$\eta\dvtx (\mathcal{X}\times\mathdsR )^{2N}\rightarrow\mathdsR $ that
is exchangeable with respect to its $2\times2N$ arguments, for any
measurable function
$\lambda\dvtx (\mathcal{X}\times\mathdsR )^{2N}\rightarrow\mathdsR _{+}^{*}$
that is exchangeable with respect to its $2\times2N$ arguments, for
any $\theta\in\varTheta$ and any $g\in\mathcal{G}$:
\begin{eqnarray*}
\lefteqn{\mathcal{P} \exp\Biggl(\frac{\lambda}{N}\sum_{i=1}^{N} \bigl\{
g \bigl[\theta(X_{i+N}),Y_{i+N} \bigr]-g \bigl[\theta(X_{i}),Y_{i} \bigr]
\bigr\}
- \frac{\lambda^{2}}{c_{g}N^{2}} \sum_{i=1}^{2N} g \bigl[\theta
(X_{i}),Y_{i} \bigr]^{2}
- \eta \Biggr) \leq\mathcal{P} \exp ( - \eta )}
\end{eqnarray*}
and the reverse inequality:
\begin{eqnarray*}
\mathcal{P} \exp\Biggl(\frac{\lambda}{N}\sum_{i=1}^{N} \bigl\{
g \bigl[\theta(X_{i}),Y_{i} \bigr]-g \bigl[\theta(X_{i+N}),Y_{i+N} \bigr]
\bigr\}
- \frac{\lambda^{2}}{c_{g}N^{2}} \sum_{i=1}^{2N} g \bigl[\theta
(X_{i}),Y_{i} \bigr]^{2}
- \eta \Biggr) \leq\mathcal{P} \exp ( - \eta ),
\end{eqnarray*}
where we write:
\begin{eqnarray*}
\lefteqn{\eta =\eta\bigl((X_{1},Y_{1}),\ldots,(X_{2N},Y_{2N}) \bigr),}
\\
\lefteqn{\lambda =\lambda\bigl((X_{1},Y_{1}),\ldots,(X_{2N},Y_{2N})
\bigr)}
\end{eqnarray*}
for short, and:
\[
c_{g}= \cases{
2 & if $g$ is nonnegative,
\cr\noalign{}
1  & otherwise.}
\]
\end{lemma}

\begin{pf}
In order to prove the first inequality, we write:
\begin{eqnarray*}
\lefteqn{\mathcal{P} \exp\Biggl(\frac{\lambda}{N}\sum_{i=1}^{N} \bigl\{
g \bigl[\theta(X_{i+N}),Y_{i+N} \bigr]-g \bigl[\theta(X_{i}),Y_{i} \bigr]
\bigr\}
- \frac{\lambda^{2}}{N^{2}} \sum_{i=1}^{2N} g \bigl[\theta(X_{i}),Y_{i} \bigr]^{2}
- \eta \Biggr)}
\\
\lefteqn{\quad = \mathcal{P} \exp\Biggl(\sum_{i=1}^{N} \log\cosh \biggl\{
\frac{\lambda}{N}
g \bigl[\theta(X_{i+N}),Y_{i+N} \bigr]-\frac{\lambda}{N}g \bigl[\theta(X_{i}),Y_{i} \bigr]
\biggr\}
- \frac{\lambda^{2}}{N^{2}} \sum_{i=1}^{2N} g \bigl[\theta(X_{i}),Y_{i} \bigr]^{2}
- \eta \Biggr).}
\end{eqnarray*}
This last step is true because $\mathcal{P}$ is exchangeable. We
conclude by using the inequality:
\[
\forall x\in\mathdsR ,\quad \log\cosh x \leq\frac{x^{2}}{2} .
\]
We obtain:
\begin{eqnarray*}
\log\cosh \biggl\{\frac{\lambda}{N}
g \bigl[\theta(X_{i+N}),Y_{i+N} \bigr]-\frac{\lambda}{N}g \bigl[\theta(X_{i}),Y_{i} \bigr]
\biggr\}
&\leq& \frac{\lambda^{2}}{2N^{2}} \bigl\{
g \bigl[\theta(X_{i+N}),Y_{i+N} \bigr]-g \bigl[\theta(X_{i}),Y_{i} \bigr] \bigr\}^{2}
\\
&\leq& \frac{\lambda^{2}}{c_{g} N^{2}} g \bigl[\theta(X_{i}),
Y_{i} \bigr]^{2} .
\end{eqnarray*}
The proof for the reverse inequality is exactly the same.
\end{pf}

We can now give the proof of the theorems.

\begin{pf*}{Proof of Theorem \ref{thtrans}}
From now on we assume that the hypothesis of Theorem \ref{thtrans}
are satisfied. Let us choose $\varepsilon'>0$ and apply Lemma
\ref{learninglemma} with $\eta=-\log\varepsilon'$, and $g$ such
that $g(u,u')=uu'$. We obtain: for any exchangeable distribution~$\mathcal{P}$, for any measurable function
$\lambda\dvt (\mathcal{X}\times\mathdsR )^{2N}\rightarrow\mathdsR _{+}^{*}$
that is exchangeable with respect to its $2\times2N$ arguments, for
any $\theta\in\varTheta$:
\[
\mathcal{P}
\exp\Biggl(\frac{\lambda}{N}\sum_{i=1}^{N} \bigl[\theta(X_{i+N})Y_{i+N}-\theta
(X_{i})Y_{i} \bigr]
- \frac{\lambda^{2}}{N^{2}} \sum_{i=1}^{2N} \theta(X_{i})^{2}Y_{i}^{2}
+ \log\varepsilon' \Biggr) \leq\varepsilon'
\]
and the reverse inequality:
\[
\mathcal{P}
\exp\Biggl(\frac{\lambda}{N}\sum_{i=1}^{N} \bigl[\theta(X_{i})Y_{i}-\theta
(X_{i+N})Y_{i+N} \bigr]
- \frac{\lambda^{2}}{N^{2}} \sum_{i=1}^{2N} \theta(X_{i})^{2}Y_{i}^{2}
+ \log\varepsilon' \Biggr) \leq\varepsilon'.
\]

Let us denote:
\[
f\bigl(\theta,\varepsilon',\lambda\bigr) = \lambda\Bigg|
\frac{1}{N}\sum_{i=1}^{N} \bigl[\theta(X_{i+N})Y_{i+N}-\theta(X_{i})Y_{i} \bigr]
- \frac{\lambda^{2}}{N^{2}} \sum_{i=1}^{2N} \theta(X_{i})^{2}Y_{i}^{2} \Bigg|
+ \log\varepsilon'.
\]

The previous inequalities imply that: for any exchangeable
$\mathcal{P}$, for any measurable function
$\lambda\dvtx (\mathcal{X}\times\mathdsR )^{2N}\rightarrow\mathdsR _{+}^{*}$
that is exchangeable with respect to its $2\times2N$ arguments, for
any $ \theta\in\varTheta$:
%
\begin{equation}\label{ineq1}
\mathcal{P} \exp
f\bigl((Z_{1},\ldots,Z_{2N}),\theta,\varepsilon',\lambda\bigr) \leq2
\varepsilon' .
\end{equation}

Now, let us introduce a new conditional probability measure:
\[
\overline{P} = \frac{1}{(2N)!} \sum_{\sigma\in\mathfrak{S}_{2N}} \delta
_{(X_{\sigma_{i}},Y_{\sigma_{i}})_{i\in\{1,\ldots,2N\}}} .
\]

Remark that $P_{2N}$ being exchangeable, we have, for any bounded
function $h\dvtx  (\mathcal{X}\times\mathdsR  )^{2N}
\rightarrow\mathdsR $,
\[
P_{2N} h = P_{2N} (\overline{P} h ) .
\]
The measure $\overline{P}$ is exchangeable, so we can apply Eq.~(\ref{ineq1}).
For any values of $Z_{1},\ldots,Z_{2N}$ we have:
\[
\forall\theta\in\varTheta,\quad \overline{P} \exp f\bigl((Z_{1},\ldots
,Z_{2N}),\theta,\varepsilon',\lambda\bigr) \leq2 \varepsilon' .
\]

In particular, we can choose $\theta=\theta(Z_{1},\ldots,Z_{2N})$ as an
exchangeable function of $(Z_{1},\ldots,Z_{2N})$, because we will have:
\begin{eqnarray*}
\lefteqn{\frac{1}{(2N)!} \sum_{\sigma\in\mathfrak{S}_{2N}} \exp
f\bigl((Z_{\sigma(1)},\ldots,Z_{\sigma(2N)}),\theta(Z_{\sigma(1)},\ldots
,Z_{\sigma(2N)}),\varepsilon',\lambda\bigr)}
\\
\lefteqn{\quad = \frac{1}{(2N)!} \sum_{\sigma\in\mathfrak{S}_{2N}} \exp
f\bigl((Z_{\sigma(1)},\ldots,Z_{\sigma(2N)}),\theta(Z_{1},\ldots
,Z_{2N}),\varepsilon',\lambda\bigr)
\leq\varepsilon' .}
\end{eqnarray*}

Here, we choose as functions $\theta$ the members of $\varTheta_{0}$:
$\theta_{1},\ldots,\theta_{m}$ (remember that we choose this indexation
in such a way that for any $k$, $ \theta_{k}$ is an exchangeable
function of $(Z_{1},\ldots,Z_{2N})$). We have, for any
$\lambda_{1},\ldots,\lambda_{m}$ that are $m$ exchangeable functions of
$(Z_{1},\ldots,Z_{2N})$:
\begin{eqnarray*}
\lefteqn{P_{2N} \bigl[\exists k\in\{1,\ldots,m\},
f\bigl((Z_{1},\ldots,Z_{2N}),\theta_{k},\varepsilon',\lambda_{k}\bigr) > 0
\bigr]}
\\
\lefteqn{\quad = P_{2N} \Biggl[\bigcup_{k=1}^{m} \bigl\{
f\bigl((Z_{1},\ldots,Z_{2N}),\theta_{k},\varepsilon',\lambda_{k}\bigr) > 0
\bigr\} \Biggr]}
\\
\lefteqn{\quad \leq P_{2N} \Biggl[\sum_{k=1}^{m} 1 \bigl(
f\bigl((Z_{1},\ldots,Z_{2N}),\theta_{k},\varepsilon',\lambda_{k}\bigr) > 0
\bigr) \Biggr]}
\\
\lefteqn{\quad = P_{2N}\overline{P} \Biggl[\sum_{k=1}^{m} 1 \bigl(
f\bigl((Z_{1},\ldots,Z_{2N}),\theta_{k},\varepsilon',\lambda_{k}\bigr) > 0
\bigr) \Biggr]}
\\
\lefteqn{\quad = P_{2N} \sum_{k=1}^{m} \overline{P} \bigl[ 1 \bigl(
f\bigl((Z_{1},\ldots,Z_{2N}),\theta_{k},\varepsilon',\lambda_{k}\bigr) > 0
\bigr) \bigr]}
\\
\lefteqn{\quad \leq P_{2N} \sum_{k=1}^{m} \overline{P} \exp
f\bigl((Z_{1},\ldots,Z_{2N}),\theta_{k},\varepsilon',\lambda_{k}\bigr).}
\end{eqnarray*}

Now let us apply inequality (\ref{ineq1}), we obtain:
\[
P_{2N} \bigl[\exists k\in\{1,\ldots,m\},
f\bigl((Z_{1},\ldots,Z_{2N}),\theta_{k},\varepsilon',\lambda_{k}\bigr) > 0
\bigr] \leq P_{2N} \sum_{k=1}^{m} 2 \varepsilon' = 2 \varepsilon' m
= \varepsilon
\]
if we choose:
\[
\varepsilon'=\frac{\varepsilon}{2m} .
\]

From now, we assume that the event:
\[
\biggl\{ \forall k\in\{1,\ldots,m\}, f \biggl((Z_{1},\ldots,Z_{2N}),\theta_{k},\frac
{\varepsilon}{2m},\lambda_{k} \biggr)\leq0 \biggr\}
\]
is satisfied. It can be written, for any $k\in\{1,\ldots,m\}$:
\[
\Bigg|\frac{1}{N}\sum_{i=1}^{N} \bigl[\theta_{k}(X_{i+N})Y_{i+N}-\theta
_{k}(X_{i})Y_{i} \bigr] \Bigg|
\leq\frac{\lambda_{k}}{N^{2}}
\sum_{i=1}^{2N}\theta_{k}(X_{i})^{2}Y_{i}^{2} +
\frac{\log ({2m}/{\varepsilon})}{\lambda_{k}} .
\]

Let us divide both inequalities by:
\[
\frac{1}{N}\sum_{i=N+1}^{2N}\theta_{k}(X_{i})^{2} .
\]
We obtain, for any $ k\in\{1,\ldots,m\}$:
\[
\big|\alpha_{2}^{k} - \mathcal{C}^{k} \alpha_{1}^{k} \big| \leq
\frac{ (\lambda_{k}/{N^{2}})
\sum_{i=1}^{2N}\theta_{k}(X_{i})^{2}Y_{i}^{2} +
({\log (2m/{\varepsilon})})/{\lambda_{k}}
}{({1}/{N})\sum_{i=N+1}^{2N}\theta_{k}(X_{i})^{2}} .
\]

It is now time to choose the functions $\lambda_{k}$. We try to
optimize the right-hand side with respect to $\lambda_{k}$, and
obtain a minimal value for:
\[
\lambda_{k} = \sqrt{\frac{N\log ({2m}/{\varepsilon})}{({1}/{N})\sum
_{i=1}^{2N}\theta_{k}(X_{i})^{2}Y_{i}^{2}}} .
\]
This choice is admissible because it is exchangeable with respect to
$(Z_{1},\ldots,Z_{2N})$.

So we have, for any $k\in\{1,\ldots,m\}$:
\[
\big| \mathcal{C}^{k} \alpha_{1}^{k} - \alpha_{2}^{k} \big|
\leq2 \frac{
\sqrt{({1}/{N^{2}})\sum_{i=1}^{2N} [\theta_{k}(X_{i})^{2}Y_{i}^{2} ]
\log ({2m}/{\varepsilon})}
}{({1}/{N})\sum_{i=N+1}^{2N}\theta_{k}(X_{i})^{2}} .
\]

Finally, remark that:
\[
\big| \mathcal{C}^{k} \alpha_{1}^{k} - \alpha_{2}^{k} \big| =
\sqrt{\frac{r_{2} [(\mathcal{C}^{k}
\alpha_{1}^{k})\theta_{k} ] - r_{2}(\alpha_{2}^{k}\theta_{k}) }
{({1}/{N})\sum_{i=N+1}^{2N}\theta_{k}(X_{i})^{2}}} ,
\]
which leads
to the conclusion that for any $ k\in\{1,\ldots,m\}$:
\[
r_{2} \bigl[\bigl(\mathcal{C}^{k} \alpha_{1}^{k}\bigr)\theta_{k} \bigr] -
r_{2}\bigl(\alpha_{2}^{k}\theta_{k}\bigr) \leq2^{2} \frac{
({1}/{N^{2}})\sum_{i=1}^{2N} [\theta_{k}(X_{i})^{2}Y_{i}^{2} ]
\log ({2m}/{\varepsilon})
}{({1}/{N})\sum_{i=N+1}^{2N}\theta_{k}(X_{i})^{2}} .
\]
This ends the proof.
\end{pf*}

\begin{pf*}{Proof of Theorem \ref{cortrans}}
We write:
\[
\frac{1}{N}\sum_{i=1}^{2N} \theta_{k}(X_{i})^{2}Y_{i}^{2}
= \frac{1}{N}\sum_{i=1}^{N} \theta_{k}(X_{i})^{2}Y_{i}^{2}
+ \frac{1}{N}\sum_{i=N+1}^{2N} \theta_{k}(X_{i})^{2}Y_{i}^{2}
\]
and try to upper bound the second term. We apply Lemma
\ref{learninglemma}, but this time with $g$ such that
$g(u)=(uu')^{2}$ that is nonnegative, and obtain, for any
$\varepsilon$, for any (exchangeables) $\theta$ and $\lambda$:
\[
\frac{1}{N}\sum_{i=N+1}^{2N} \theta_{k}(X_{i})^{2}Y_{i}^{2} \leq
\frac{1}{N}\sum_{i=1}^{N} \theta_{k}(X_{i})^{2}Y_{i}^{2} +
\frac{\lambda}{2N}\biggl(\frac{1}{N}\biggr)\sum_{i=1}^{2N}
\theta_{k}(X_{i})^{4}Y_{i}^{4} + \frac{\log\varepsilon}{\lambda} .
\]
We choose:
\[
\lambda= \sqrt{\frac{2N \log\varepsilon}{({1}/{N})\sum_{i=1}^{2N}
\theta_{k}(X_{i})^{4}Y_{i}^{4}}} ,
\]
we apply this result to every $\theta\in\varTheta_{0}$, and combine it
with Theorem \ref{thtrans} by a union bound argument to obtain the
result.
\end{pf*}

\subsection[Proof of Theorem 3.3]{Proof of Theorem \textup{{\protect\ref{thimp1}}}}
\label{subproof3}

First of all, we give the following obvious variant of Lemma~\ref{learninglemma}:
\begin{lemma} \label{ll2}
For any exchangeable probability distribution $\mathcal{P}$ on
$ (Z_{1},\ldots,Z_{2N} )$, for any measurable function
$\eta\dvtx (\mathcal{X}\times\mathdsR )^{2N}\rightarrow\mathdsR $ that
is exchangeable with respect to its $2\times2N$ arguments, for any
measurable function
$\lambda\dvtx (\mathcal{X}\times\mathdsR )^{2N}\rightarrow\mathdsR _{+}^{*}$
that is exchangeable with respect to its $2\times2N$ arguments, for
any $\theta\in\varTheta$:
\begin{eqnarray*}
\lefteqn{\mathcal{P} \exp\Biggl(\frac{\lambda}{N}\sum_{i=1}^{N}
\bigl\{ \bigl[\theta(X_{i+N})Y_{i+N}-\alpha(\theta)\theta(X_{i+N})^{2} \bigr]
- \bigl[\theta(X_{i})Y_{i}-\alpha(\theta)\theta(X_{i})^{2} \bigr]
\bigr\}}
\\
\lefteqn{\quad{} - \frac{\lambda^{2}}{N^{2}} \sum_{i=1}^{2N} \bigl[\theta(X_{i})Y_{i}-\alpha
(\theta)\theta(X_{i})^{2} \bigr]^{2}- \eta \Biggr) \leq\mathcal{P} \exp ( - \eta )}
\end{eqnarray*}
and the reverse inequality, where:
\[
\alpha(\theta) = \mathop{\arg\min}_{\alpha\in\mathdsR } r_{1,2} (\alpha\theta) .
\]
\end{lemma}

\begin{pf}
This is actually just an application of Lemma~\ref{learninglemma},
we just need to remark that $\alpha(\theta)$ is an exchangeable
function of $(Z_{1},\ldots,Z_{2N})$, and so we can take in Lemma~\ref{learninglemma}:
\[
g\bigl(u,u'\bigr) = u u' - u^{2} \alpha(\theta) ,
\]
that means that:
\[
g \bigl[\theta(X_{i}),Y_{i} \bigr] = \theta(X_{i})Y_{i}-\alpha(\theta)\theta
(X_{i})^{2} .
\]\upqed
\end{pf}

\begin{pf*}{Proof of Theorem~\ref{thimp1}}
Proceeding exactly in the same way as in the proof of Theorem
\ref{thtrans}, we obtain the following inequality with probability
at least $1-\varepsilon$:
%
\begin{equation} \label{var1}
r_{2}\bigl(\mathcal{C}^{k}\alpha_{1}^{k}\theta_{k}\bigr)-r_{2}\bigl(\alpha
_{2}^{k}\theta_{k}\bigr)
\leq
4 \biggl[\frac{({1}/{N})\sum_{i=1}^{2N} [\theta_{k}(X_{i})Y_{i}-\alpha
_{1,2}^{k}\theta_{k}(X_{i})^{2} ]^{2}}
{({1}/{N})\sum_{i=N+1}^{2N}\theta_{k}(X_{i})^{2}} \biggr]\frac{\log (2m/\varepsilon)}{N}.
\end{equation}
This proves the theorem.
\end{pf*}

Before giving the proof of the next theorem, let us see how we can
make the first-order term observable in this theorem. For example,
we can write:
\begin{eqnarray*}
\bigl[\theta_{k}(X_{i})Y_{i}-\alpha_{1,2}^{k}\theta_{k}(X_{i})^{2} \bigr]^{2}
&=& \bigl[\theta_{k}(X_{i})Y_{i}-\alpha_{1}^{k}\theta_{k}(X_{i})^{2} \bigr]^{2}
+ \bigl[\alpha_{1}^{k}-\alpha_{1,2}^{k} \bigr]^{2}\theta_{k}(X_{i})^{4}
\\
&&{}+2 \bigl[\theta_{k}(X_{i})Y_{i}-\alpha_{1}^{k}\theta_{k}(X_{i})^{2} \bigr]
\bigl[\alpha_{1}^{k}-\alpha_{1,2}^{k} \bigr]\theta_{k}(X_{i})^{2}.
\end{eqnarray*}

Remark that it is obvious that:
\[
\big| \alpha_{1}^{k}-\alpha_{1,2}^{k} \big| \leq\big| \alpha_{1}^{k}-\alpha
_{2}^{k} \big| ,
\]
and so:
\begin{eqnarray*}
\bigl[\theta_{k}(X_{i})Y_{i}-\alpha_{1,2}^{k}\theta_{k}(X_{i})^{2} \bigr]^{2}
&\leq& \bigl[\theta_{k}(X_{i})Y_{i}-\alpha_{1}^{k}\theta_{k}(X_{i})^{2} \bigr]^{2}
+ \bigl[\alpha_{1}^{k}-\alpha_{2}^{k} \bigr]^{2}\theta_{k}(X_{i})^{4}
\\
&&{} + 2 \big|\theta_{k}(X_{i})Y_{i}-\alpha_{1}^{k}\theta_{k}(X_{i})^{2} \big|
\big|\alpha_{1}^{k}-\alpha_{2}^{k} \big|\theta_{k}(X_{i})^{2}.
\end{eqnarray*}

Now, just write:
\[
\alpha_{1}^{k}-\alpha_{2}^{k} = \bigl(1-\mathcal{C}^{k} \bigr)\alpha_{1}^{k} -
\bigl(\mathcal{C}^{k}\alpha_{1}^{k}-\alpha_{2}^{k} \bigr)
\]
and so we get:
\begin{eqnarray*}
\bigl[\theta_{k}(X_{i})Y_{i}-\alpha_{1,2}^{k}\theta_{k}(X_{i})^{2} \bigr]^{2}
&\leq&
\bigl[\theta_{k}(X_{i})Y_{i}-\alpha_{1}^{k}\theta_{k}(X_{i})^{2} \bigr]^{2}
+ \bigl[\mathcal{C}^{k}\alpha_{1}^{k}-\alpha_{2}^{k} \bigr]^{2}\theta_{k}(X_{i})^{4}
\\
&&{} + 2 \big|\mathcal{C}^{k}\alpha_{1}^{k}-\alpha_{2}^{k} \big| \big|\bigl(1-\mathcal
{C}^{k}\bigr)\alpha_{1}^{k} \big|\theta_{k}(X_{i})^{4}
+ \bigl(1-\mathcal{C}^{k} \bigr)^{2} \bigl(\alpha_{1}^{k} \bigr)^{2}
\theta_{k}(X_{i})^{4}
\\
&&{} + 2 \big|\theta_{k}(X_{i})Y_{i}-\alpha_{1}^{k}\theta_{k}(X_{i})^{2} \big|
\big|\mathcal{C}^{k}\alpha_{1}^{k}-\alpha_{2}^{k} \big|\theta_{k}(X_{i})^{2}
\\
&&{} + 2 \big|\theta_{k}(X_{i})Y_{i}-\alpha_{1}^{k}\theta_{k}(X_{i})^{2} \big|
\big|\bigl(\mathcal{C}^{k}-1\bigr)\alpha_{1}^{k} \big|\theta_{k}(X_{i})^{2}.
\end{eqnarray*}
So finally, Eq.~(\ref{var1}) left us with a second degree
inequality with respect to
$ |\mathcal{C}^{k}\alpha_{1}^{k}-\alpha_{2}^{k} |$ or
$r_{2}(\mathcal{C}^{k}\alpha_{1}^{k}\theta_{k})-r_{2}(\alpha
_{2}^{k}\theta_{k})$ that we can solve to obtain the following result: with probability
at least $1-\varepsilon$, as soon as we have:
\[
\Biggl[\frac{1}{N}\sum_{i=N+1}^{2N}\theta_{k}(X_{i})^{2} \Biggr]^{2}
> \Biggl[\frac{1}{N}\sum_{i=1}^{2N}\theta_{k}(X_{i})^{4} \Biggr] \frac{4\log
({2m}/{\varepsilon})}{N} ,
\]
which is always true for large enough $N$, the quantity
$ |\mathcal{C}^{k}\alpha_{1}^{k}-\alpha_{2}^{k} |$ belongs
to the interval:
\[
\Biggl[
\frac{2\log (2m/{\varepsilon})}{N} \frac{b \pm\sqrt{b^{2} +
a ( (N/{\log (2m/\varepsilon)})
[(1/N)\sum_{i=N+1}^{2N}\theta_{k}(X_{i})^{2} ]^{2}
-(4/N)\sum_{i=1}^{2N}\theta_{k}(X_{i})^{4} )}
}{ [(1/N)\sum_{i=N+1}^{2N}\theta_{k}(X_{i})^{2} ]^{2}
- (4\log (2m/\varepsilon)/N)
[(1/N)\sum_{i=1}^{2N}\theta_{k}(X_{i})^{4} ] }
\Biggr]
\]
with the following notations:
\begin{eqnarray*}
\lefteqn{a  = \frac{1}{N}
\sum_{i=1}^{2N} \bigl[ \big|\theta_{k}(X_{i})Y_{i}-\alpha_{1}^{k}\theta
_{k}(X_{i})^{2} \big|
+ \big|\alpha_{k}^{1}\bigl(1-\mathcal{C}^{k}\bigr) \big|\theta_{k}(X_{i})^{2}
\bigr]^{2},}
\\
\lefteqn{b  = \frac{1}{N}\sum_{i=1}^{2N} 2 \theta_{k}(X_{i})^{2}
\bigl[ \big|\alpha_{k}^{1}\bigl(1-\mathcal{C}^{k}\bigr) \big|\theta_{k}(X_{i})^{2}
+ \big|\theta_{k}(X_{i})Y_{i}-\alpha_{1}^{k}\theta_{k}(X_{i})^{2} \big|
\bigr].}
\end{eqnarray*}
Remark that only one of the bounds of the interval is positive. So
we obtain the following result: with $P_{2N}$-probability at least
$1-\varepsilon$, as soon as:
\[
\Biggl[\frac{1}{N}\sum_{i=N+1}^{2N}\theta_{k}(X_{i})^{2} \Biggr]^{2}
> \Biggl[\frac{1}{N}\sum_{i=1}^{2N}\theta_{k}(X_{i})^{4} \Biggr] \frac{4\log (2m/\varepsilon)}{N}
\]
we have:
\begin{eqnarray*}
\lefteqn{\forall k\in\{1,\ldots,m\},}
\\
\lefteqn{ r_{2} \bigl[\bigl(\mathcal{C}^{k} \alpha_{1}^{k}\bigr)\theta_{k} \bigr]-r_{2}\bigl(\alpha_{2}^{k}\theta_{k}\bigr)}
\\
\lefteqn{\quad \leq
\frac{4\log^{2} ({2m}/{\varepsilon})}{N^{2}} \Biggl[\frac{1}{N}\sum
_{i=1}^{2N}\theta_{k}(X_{i})^{2} \Biggr]}
\\
\lefteqn{\qquad {}\times \Biggl[ \frac{ b + \sqrt{b^{2} + a (
(N/\log(2m/\varepsilon))
[({1}/{N})\sum_{i=N+1}^{2N}\theta_{k}(X_{i})^{2} ]^{2}
-({4}/{N})\sum_{i=1}^{2N}\theta_{k}(X_{i})^{4} )}
}{ [({1}/{N})\sum_{i=N+1}^{2N}\theta_{k}(X_{i})^{2} ]^{2}
- ((4\log (2m/\varepsilon))/N)
[(1/N)\sum_{i=1}^{2N}\theta_{k}(X_{i})^{4} ] }
\Biggr]^{2}.}
\end{eqnarray*}
We can notice that this bound may be written:
\begin{eqnarray*}
r_{2} \bigl[\bigl(\mathcal{C}^{k}
\alpha_{1}^{k}\bigr)\theta_{k} \bigr]-r_{2}\bigl(\alpha_{2}^{k}\theta_{k}\bigr)
&\leq& \frac{8 a \log ({2m}/{\varepsilon})}{N} +
\EuScript{O} \biggl( \biggl[\frac{\log ({m}/{\varepsilon})}{N} \biggr]^{{3}/{2}} \biggr)
\\
&=&  \frac{8\log (2m/\varepsilon)}{N}
\Biggl[\frac{1}{N}\sum_{i=1}^{2N} \bigl(\theta_{k}(X_{i})Y_{i}-\alpha
_{1}^{k}\theta_{k}(X_{i})^{2} \bigr)^{2} \Biggr]
+\EuScript{O} \biggl( \biggl[\frac{\log ({m}/{\varepsilon})}{N} \biggr]^{{3}/{2}} \biggr).
\end{eqnarray*}
The next step would be now to replace the bound by an observable
quantity, by getting a bound like:
\[
\frac{1}{N}\sum_{i=1}^{2N} \bigl(\theta_{k}(X_{i})Y_{i}-\alpha_{1}^{k}\theta
_{k}(X_{i})^{2} \bigr)^{2}
\leq\frac{2}{N}\sum_{i=1}^{N} \bigl(\theta_{k}(X_{i})Y_{i}-\alpha
_{1}^{k}\theta_{k}(X_{i})^{2} \bigr)^{2} +
\EuScript{O} \biggl(\frac{\log ({m}/{\varepsilon})}{N} \biggr)
\]
with high probability. This can be done very simply, using Lemma
\ref{learninglemma} with this time:
\[
g\bigl(u,u'\bigr) = \bigl(uu'-u^{2}\alpha(\theta) \bigr)^{2} .
\]
We obtain the bound:
\[
r_{2} \bigl[\bigl(\mathcal{C}^{k}
\alpha_{1}^{k}\bigr)\theta_{k} \bigr]-r_{2}\bigl(\alpha_{2}^{k}\theta_{k}\bigr)
\leq\frac{16 \log (4m/\varepsilon)}{N}
\Biggl[\frac{1}{N}\sum_{i=1}^{N} \bigl(\theta_{k}(X_{i})Y_{i}-\alpha_{1}^{k}\theta
_{k}(X_{i})^{2} \bigr)^{2} \Biggr]
+\EuScript{O} \biggl( \biggl[\frac{\log (m/\varepsilon)}{N} \biggr]^{{3}/{2}} \biggr).
\]

\subsection[Proof of Theorem 3.4]{Proof of Theorem \textup{{\protect\ref{thimp2}}}}
\label{subproof4}

The proof is exactly similar, we just use a new variant of lemma
\ref{learninglemma}, that is based on an idea introduced by Catoni
\cite{CatVapnik} in the context of classification.

\begin{dfn}
Let us write:
\[
T_{\theta}(Z_{i}) = \theta(X_{i}) Y_{i}
\]
for short. We also introduce a conditional probability measure:
\[
\mathcal{P}^{(2)}=\frac{1}{N!}\sum_{\sigma\in\mathfrak{S}_{N}}\delta
_{(Z_{1},\ldots,Z_{N},Z_{N+\sigma(1)},\ldots,Z_{N+\sigma(N)})} .
\]
\end{dfn}

Remark that, because $\mathcal{P}$ is exchangeable, we have, for any
function $h$:
\[
\mathcal{P}h=\mathcal{P} \bigl[\mathcal{P}^{(2)}h \bigr] .
\]

\begin{lemma} \label{lemma2}
For any exchangeable probability distribution $\mathcal{P}$ on
$ (Z_{1},\ldots,Z_{2N} )$, for any measurable function
$\eta\dvtx (\mathcal{X}\times\mathdsR )^{2N}\rightarrow\mathdsR $ that
is exchangeable with respect to its $2\times2N$ arguments, for any
measurable function
$\lambda\dvtx (\mathcal{X}\times\mathdsR )^{2N}\rightarrow\mathdsR _{+}^{*}$
which is such that, for any $i\in\{1,\ldots,2N\}$:
\[
\lambda(Z_{1},\ldots,Z_{2N}) = \lambda(Z_{1},\ldots
,Z_{i-1},Z_{i+N},Z_{i+1},\ldots,Z_{i+N-1},Z_{i},Z_{i+N+1},\ldots,Z_{2N}),
\]
for any $\theta\in\varTheta$:
\[
\mathcal{P}
\exp\Biggl\{\frac{\mathcal{P}^{(2)}\lambda}{N}\sum_{i=1}^{N} \bigl[T_{\theta
}(Z_{i})-T_{\theta}(Z_{i+N}) \bigr]
- \mathcal{P}^{(2)} \Biggl[\frac{\lambda^{2}}{2N^{2}} \frac{1}{N}\sum
_{i=1}^{N} \bigl[T_{\theta}(Z_{i})-T_{\theta}(Z_{i+N}) \bigr]^{2} \Biggr]
- \eta \Biggr\} \leq\mathcal{P} \exp ( - \eta )
\]
and the reverse inequality.
\end{lemma}

\begin{pf}
Let $\mathcal{L}hs$ denote the left-hand side of Lemma \ref{lemma2}.
For short, let us put:
\[
s(\theta) = \frac{1}{N}\sum_{i=1}^{N} \bigl[\theta(X_{i+N})Y_{i+N}-\theta
(X_{i})Y_{i} \bigr]^{2}
= \frac{1}{N}\sum_{i=1}^{N} \bigl[T_{\theta}(Z_{i})-T_{\theta}(Z_{i+N}) \bigr]^{2}.
\]
Then we have:
\begin{eqnarray*}
\mathcal{L}hs  &=&  P_{2N}\exp
P^{(2)} \Biggl(\frac{\lambda}{N}\sum_{i=1}^{N} \bigl[T_{\theta}(Z_{i})-T_{\theta
}(Z_{i+N}) \bigr]
-\frac{\lambda^{2}}{2N}s(\theta) -\eta\Biggr)
\\
&\leq&  P_{2N} P^{(2)}\exp \Biggl(\frac{\lambda}{N}\sum_{i=1}^{N}
\bigl[T_{\theta}(Z_{i})-T_{\theta}(Z_{i+N}) \bigr]
-\frac{\lambda^{2}}{2N}s(\theta) -\eta\Biggr),
\end{eqnarray*}
by Jensen's conditional inequality. Now, we can conclude as in Lemma
\ref{learninglemma}:
\begin{eqnarray*}
\mathcal{L}hs  &=&
P_{2N}\exp\Biggl(\sum_{i=1}^{N}\log\cosh\biggl\{\frac{\lambda}{N} \bigl[T_{\theta
}(Z_{i})-T_{\theta}(Z_{i+N}) \bigr] \biggr\}
-\frac{\lambda^{2}}{2N}s(\theta)-\eta\Biggr)
\\
& \leq& P_{2N}\exp\Biggl(\frac{\lambda^{2}}{2N^{2}}\sum_{i=1}^{N}
\bigl[T_{\theta}(Z_{i})-T_{\theta}(Z_{i+N}) \bigr]^{2}
-\frac{\lambda^{2}}{2N}s(\theta)-\eta\Biggr)
\\
&=&  P_{2N}\exp (-\eta).
\end{eqnarray*}\upqed
\end{pf}

\begin{pf*}{Proof of Theorem \ref{thimp2}}
We apply both inequalities of Lemma \ref{lemma2} to every
$\theta_{k},k\in\{1,\ldots,m\}$, and we take:
\[
\lambda= \sqrt{\frac{2N\log (2m/\varepsilon)}{s(\theta)}} .
\]
We obtain, for any $k\in\{1,\ldots,m\}$:
\[
\mathcal{P}
\exp\Biggl\{\frac{\mathcal{P}^{(2)}\lambda}{N}\sum_{i=1}^{N} \bigl[T_{\theta
}(Z_{i})-T_{\theta}(Z_{i+N}) \bigr]
- \log\frac{2m}{\varepsilon}
- \eta \Biggr\} \leq\varepsilon.
\]
Or, with probability at least $1-\varepsilon$, for any $k$:
\[
\frac{1}{N}\sum_{i=1}^{N} \bigl[T_{\theta}(Z_{i})-T_{\theta}(Z_{i+N}) \bigr]
\leq\sqrt{\frac{2\log (2m/\varepsilon)}{N}}
\bigl[\mathcal{P}^{(2)} \bigl(s(\theta)^{-{1}/{2}} \bigr) \bigr]^{-1},
\]
so:
\[
\Biggl[\frac{1}{N}\sum_{i=1}^{N}T_{\theta}(Z_{i})-\frac{1}{N}\sum
_{i=N+1}^{2N}T_{\theta}(Z_{i}) \Biggr]^{2}
\leq\frac{2\log ({2m}/{\varepsilon})}{N}
\mathcal{P}^{(2)}s(\theta).
\]
We end the first part of the proof by noting that:
\[
\mathcal{P}^{(2)}s(\theta)
=
V_{1}(\theta)+V_{2}(\theta)+ \Biggl[\frac{1}{N}\sum_{i=1}^{N}T_{\theta
}(Z_{i})-\frac{1}{N}\sum_{i=N+1}^{2N}T_{\theta}(Z_{i}) \Biggr]^{2}.
\]
Now, let us see how we can obtain the second part of the theorem.
Note that:
\[
V_{2}(\theta)=\frac{1}{N}\sum_{i=N+1}^{2N}T_{\theta}(Z_{i})^{2} -
\Biggl(\frac{1}{N}\sum_{i=N+1}^{2N}T_{\theta}(Z_{i}) \Biggr)^{2}.
\]
We upper bound the first term by using Lemma \ref{learninglemma}
with
$g(\theta(X_{i}),Y_{i})=\theta(X_{i})^{2}Y_{i}^{2}=T_{\theta}(Z_{i})^{2}$,
so with probability at least $1-\varepsilon$, for any $k$:
\[
\frac{1}{N}\sum_{i=N+1}^{2N}T_{\theta}(Z_{i})^{2} \leq
\frac{1}{N}\sum_{i=1}^{N}T_{\theta}(Z_{i})^{2} + \sqrt{\frac{2\log ({m}/{\varepsilon})
({1}/{N})\sum_{i=1}^{2N}T_{\theta}(Z_{i})^{4}}{N}}.
\]
For the second-order term, we use both inequalities of Lemma
\ref{learninglemma} with
$g(\theta(X_{i}),Y_{i})=\theta(X_{i})Y_{i}=T_{\theta}(Z_{i})$, so
with probability at least $1-\varepsilon$, for any $k$:
\begin{eqnarray*}
\Biggl(\frac{1}{N}\sum_{i=1}^{N}T_{\theta}(Z_{i}) \Biggr)^{2} -
\Biggl(\frac{1}{N}\sum_{i=N+1}^{2N}T_{\theta}(Z_{i}) \Biggr)^{2}
&\leq& \Bigg|\frac{1}{N}\sum_{i=1}^{N}T_{\theta}(Z_{i}) -
\frac{1}{N}\sum_{i=N+1}^{2N}T_{\theta}(Z_{i}) \Bigg|
\Bigg|\frac{1}{N}\sum_{i=1}^{2N}T_{\theta}(Z_{i}) \Bigg|
\\
&\leq&
2\sqrt{\frac{({1}/{N})\sum_{i=1}^{2N}T_{\theta}(Z_{i})^{2}\log
({2m}/{\varepsilon})}{N}}
\frac{1}{N}\sum_{i=1}^{2N} \big|T_{\theta}(Z_{i}) \big|.
\end{eqnarray*}
Putting all pieces together (and replacing $\varepsilon$ by
$\varepsilon/3$) ends the proof.
\end{pf*}

\subsection[Proof of Theorem 3.5]{Proof of Theorem \textup{{\protect\ref{thimp3}}}}
\label{subproof5}
\begin{pf*}{Proof of Theorem \ref{thimp3}}
We introduce the following conditional probability measures, for any
$i\in\{1,\ldots,N\}$:
{\mathindent=0pt
\begin{eqnarray*}
\PP_{i}
&=&  \frac{1}{(k+1)!}
\\
&&\!{}\times \!\!\! \sum_{\sigma\in\mathfrak{S}_{k+1}}\!
\delta_{(Z_{1},\ldots,Z_{i-1},Z_{N(\sigma(1)-1)+i},Z_{i+1},\ldots,
Z_{N+i-1},
Z_{N(\sigma(2)-1)+i},Z_{N+i+1},\ldots,\ldots,Z_{kN+i-1},Z_{N(\sigma
(k+1)-1)+i},Z_{kN+i+1},\ldots,Z_{(k+1)N})}.
\end{eqnarray*}}%
and
\[
\PP = \bigotimes_{i=1}^{N} \PP_{i}
\]
and, finally, remember that:
\[
\mathbf{P} = \frac{1}{(k+1)N}\sum_{i=1}^{(k+1)N} \delta_{Z_{i}} .
\]

Note that, by exchangeability, for any nonnegative function
\[
h\dvtx (\mathcal{X}\times\mathdsR )^{(k+1)N}\rightarrow\mathdsR
\]
we have, for any $i\in\{1,\ldots,N\}$:
\[
P_{(k+1)N} \PP_{i} h(Z_{1},\ldots,Z_{2N}) = P_{(k+1)N} h
(Z_{1},\ldots,Z_{2N}) .
\]

\begin{lemma} \label{klemma}
Let $\chi$ be a function $\mathdsR \rightarrow\mathdsR $. For any
exchangeable functions $\lambda$, $\eta\dvtx
(\mathcal{X}\times\mathdsR )^{(k+1)N}\rightarrow\mathdsR _{+} $
and $\theta\dvtx (\mathcal{X}\times\mathdsR )^{(k+1)N}\rightarrow\varTheta$
we have:
\begin{eqnarray*}
\lefteqn{\PP \exp\Biggl\{\lambda\Biggl[\frac{1}{kN}\sum_{i=N+1}^{(k+1)N}
\chi\bigl[\theta(X_{i})Y_{i} \bigr]
-\frac{1}{N}\sum_{i=1}^{N} \chi \bigl[\theta(X_{i})Y_{i} \bigr]
\Biggr]-\eta\Biggr\}}
\\
\lefteqn{\quad \leq\exp(-\eta) \exp\biggl\{
\frac{\lambda^{2}(1+k)^{2}}{2Nk^{2}} \mathbf{P}
\bigl\{ \bigl[\chi\bigl(\theta(X)Y \bigr)
- \mathbf{P} \chi\bigl(\theta(X)Y \bigr) \bigr]^{2} \bigr\}}
\\
\lefteqn{\qquad{} + \frac{\lambda^{3}(1+k)^{3}}{6N^{2} k^{3}}
\Bigl[\sup_{i\in\{1,\ldots,(k+1)N\}}\chi \bigl(\theta(X_{i})Y_{i} \bigr)-\inf_{i\in\{
1,\ldots,(k+1)N\}}
\chi\bigl(\theta(X_{i})Y_{i} \bigr) \Bigr]^{3}
\biggr\} ,}
\end{eqnarray*}
where we put $ \lambda= \lambda(Z_{1},\ldots,Z_{(k+1)N})$,
$\theta=\theta(Z_{1},\ldots,Z_{(k+1)N})$ and $ \eta=
\eta(Z_{1},\ldots,Z_{(k+1)N})$ for short. We have the reverse
inequality as well.
\end{lemma}

Before giving the proof, let us introduce the following useful
notations.
\begin{dfn}
We put, for any $\theta\in\varTheta$, for any function $\chi$:
\[
\chi_{i}^{\theta} = \chi \bigl( Y_{i}\theta(X_{i}) \bigr) ,
\]
and
\[
\chi^{\theta} = \chi \bigl( Y\theta(X) \bigr)
\]
that means that:
\[
\mathbf{P} \chi^{\theta} = \frac{1}{(k+1)N}\sum_{i=1}^{(k+1)N}\chi
_{i}^{\theta} .
\]
We also put:
\[
\mathcal{S}_{\chi}(\theta) = \sup_{i\in\{1,\ldots,(k+1)N\}}\chi
_{i}^{\theta}-\inf_{i\in\{1,\ldots,(k+1)N\}}\chi_{i}^{\theta}.
\]
\end{dfn}

\begin{pf*}{Proof of the Lemma \ref{klemma}}
Remark that, for any exchangeable functions $\lambda$, $\eta\dvtx
(\mathcal{X}\times\mathdsR )^{(k+1)N}\rightarrow\mathdsR _{+} $
and $\theta\dvtx (\mathcal{X}\times\mathdsR )^{kN}\rightarrow\varTheta$ we
have:
\begin{eqnarray*}
\lefteqn{\PP \exp\Biggl\{\lambda\Biggl[\frac{1}{kN}\sum_{i=N+1}^{(k+1)N}
g \bigl[\theta(X_{i})Y_{i} \bigr]
-\frac{1}{N}\sum_{i=1}^{N} g \bigl[\theta(X_{i})Y_{i} \bigr]
\Biggr]-\eta\Biggr\}}
\\
\lefteqn{\quad =  \exp(-\eta) \prod_{i=1}^{N} \PP_{i} \exp \Biggl\{
\frac{\lambda}{kN}\sum_{j=1}^{k}\chi_{i+jN}^{\theta}
- \frac{\lambda}{N}\chi_{i}^{\theta} \Biggr\}}
\\
\lefteqn{\quad =  \exp(-\eta) \prod_{i=1}^{N} \exp \Biggl\{
\frac{\lambda}{kN}\sum_{j=0}^{k}\chi_{i+jN}^{\theta} \Biggr\}
\prod_{i=1}^{N} \PP_{i} \exp \biggl\{ - \frac{\lambda(1+k)}{kN}\chi
_{i}^{\theta} \biggr\},}
\end{eqnarray*}
where we put $ \lambda= \lambda(Z_{1},\ldots,Z_{kN})$,
$\theta=\theta(Z_{1},\ldots,Z_{kN})$ and $ \eta=
\eta(Z_{1},\ldots,Z_{kN})$ for short.

Now, we have:
\[
\log\prod_{i=1}^{N} \PP_{i} \exp \biggl\{ -
\frac{\lambda(1+k)}{kN}\chi_{i}^{\theta} \biggr\} = \sum_{i=1}^{N}
\log\PP_{i} \exp \biggl\{ -
\frac{\lambda(1+k)}{kN}\chi_{i}^{\theta} \biggr\},
\]
and, for any $i\in\{1,\ldots,N\}$:
\begin{eqnarray*}
\lefteqn{\log\PP_{i} \exp \biggl\{ -
\frac{\lambda(1+k)}{Nk}\chi_{i}^{\theta} \biggr\}}
\\
\lefteqn{\quad = - \frac{\lambda(1+k)}{Nk}\PP_{i}\chi_{i}^{\theta} +
\frac{\lambda^{2}(1+k)^{2}}{2N^{2}k^{2}}
\PP_{i} \bigl[ \bigl(\chi_{i}^{\theta}-\PP_{i}\chi_{i}^{\theta}
\bigr)^{2} \bigr]}
\\
\lefteqn{\qquad{}  - \int_{0}^{{\lambda(1+k)}/{(Nk)}} \frac{1}{2}
\biggl(\frac{\lambda(1+k)}{Nk}-\beta\biggr)^{2}
\frac{1}{\PP_{i} \exp[-\beta\chi_{i}^{\theta} ]}
\PP_{i} \biggl[ \biggl(\chi_{i}^{\theta}-\frac{\PP_{i} \{ \chi
_{i}^{\theta} \exp[-\beta\chi_{i}^{\theta} ] \}}{
\PP_{i} \exp[-\beta\chi_{i}^{\theta} ]} \biggr)^{3}\exp\bigl(-\beta\chi
_{i}^{\theta} \bigr) \biggr]\,\mathrm{d}\beta.}
\end{eqnarray*}
Note that, for any $\beta\geq0$:
\[
\frac{1}{\PP_{i} \exp[-\beta\chi_{i}^{\theta} ]}
\PP_{i} \biggl[ \biggl(\chi_{i}^{\theta}-\frac{\PP_{i} \{ \chi
_{i}^{\theta} \exp[-\beta\chi_{i}^{\theta} ] \}}{
\PP_{i} \exp[-\beta\chi_{i}^{\theta} ]} \biggr)^{3}\exp\bigl(-\beta\chi
_{i}^{\theta} \bigr) \biggr]
\leq
\Bigl[\sup_{j\in\{1,\ldots,k\}}\chi_{i+(j-1)N}^{\theta}-\inf_{j\in\{1,\ldots
,k\}}\chi_{i+(j-1)N}^{\theta} \Bigr]^{3},
\]
and so:
\begin{eqnarray*}
\log\prod_{i=1}^{N} \PP_{i} \exp \biggl\{ -
\frac{\lambda(1+k)}{Nk}\chi_{i}^{\theta} \biggr\} &\leq& -
\frac{1}{N} \sum_{i=1}^{N}
\frac{\lambda(1+k)}{k}\PP_{i}\chi_{i}^{\theta}
+ \frac{1}{N} \sum_{i=1}^{N} \frac{\lambda^{2}(1+k)^{2}}{2N k^{2}}
\PP_{i} \bigl[ \bigl(\chi_{i}^{\theta}-\PP_{i}\chi_{i}^{\theta}
\bigr)^{2} \bigr]
\\
&&{} + \frac{\lambda^{3}(1+k)^{3}}{6N^{2} k^{3}}
\Bigl[\sup_{i\in\{1,\ldots,(k+1)N\}}\chi_{i}^{\theta}-\inf_{i\in\{1,\ldots
,(k+1)N\}}\chi_{i}^{\theta} \Bigr]^{3}.
\end{eqnarray*}
Note that:
\[
\PP_{i} \chi_{i}^{\theta} = \frac{1}{k+1}\sum_{j=0}^{k}\chi
_{i+jN}^{\theta}
\]
and so:
\[
\frac{1}{N} \sum_{i=1}^{N} \PP_{i} \chi_{i}^{\theta}
= \frac{1}{(k+1)N}\sum_{i=1}^{(k+1)N}\chi_{i}^{\theta} = \mathbf{P}
\chi^{\theta} ;
\]
remark also that:
\[
\frac{1}{N} \sum_{i=1}^{N} \PP_{i} \bigl[ \bigl(\chi_{i}^{\theta}-\mathdsP _{i}\chi_{i}^{\theta} \bigr)^{2} \bigr]
\leq \frac{1}{(k+1)N} \sum_{i=1}^{(k+1)N}
\Biggl[\chi_{i}^{\theta}- \Biggl(\frac{1}{(k+1)N}\sum_{j=1}^{(k+1)N}\chi
_{j}^{\theta} \Biggr) \Biggr]^{2}
= \mathbf{P} \bigl[ \bigl(\chi^{\theta}- \mathbf{P}
\chi^{\theta} \bigr)^{2} \bigr] ,
\]
we obtain:
\begin{eqnarray*}
\lefteqn{\PP \exp\Biggl\{\lambda\Biggl[\frac{1}{kN}\sum_{i=N+1}^{(k+1)N}
\theta(X_{i})Y_{i}
-\frac{1}{N}\sum_{i=1}^{N} \theta(X_{i})Y_{i} \Biggr]-\eta\Biggr\}}
\\
\lefteqn{\quad = \exp(-\eta) \exp\biggl\{
\frac{\lambda^{2}(1+k)^{2}}{2Nk^{2}} \mathbf{P}
\bigl[ \bigl(\chi^{\theta}
- \mathbf{P} \chi^{\theta} \bigr)^{2} \bigr]
+ \frac{\lambda^{3}(1+k)^{3}}{6N^{2} k^{3}}
\Bigl[\sup_{i\in\{1,\ldots,(k+1)N\}}\chi_{i}^{\theta}-\inf_{i\in\{1,\ldots
,(k+1)N\}}\chi_{i}^{\theta} \Bigr]^{3}
\biggr\} .}
\end{eqnarray*}
The proof of the reverse inequality is exactly the same.
\end{pf*}

Let us choose here again $\chi$ such that $\chi(u)=u$, namely:
$\chi=id$. By the use of a union bound argument on elements of
$\varTheta_{0}$ we obtain, for any $\varepsilon>0$, for any
exchangeable function $\lambda\dvtx
(\mathcal{X}\times\mathdsR )^{(k+1)N}\rightarrow\mathdsR _{+} $,
with probability at least $1-\varepsilon$, for any
$h\in\{1,\ldots,m\}$:
\begin{eqnarray*}
\lefteqn{\frac{1}{kN}\sum_{i=N+1}^{(k+1)N} \theta_{h}(X_{i})Y_{i}
-\frac{1}{N}\sum_{i=1}^{N} \theta_{h}(X_{i})Y_{i}}
\\
\lefteqn{\quad \leq\frac{\lambda (1+ {1}/{k} )^{2}}{2N} \mathbf{P}
\bigl[ \bigl(\chi^{\theta_{h}}
- \mathbf{P} \chi^{\theta_{h}} \bigr)^{2} \bigr]
+ \frac{\lambda^{2} (1+ {1}/{k} )^{3}}{6N^{2}}
\mathcal{S}_{id}(\theta_{h})^{3} +
\frac{\log ({m}/{\varepsilon})}{\lambda} .}
\end{eqnarray*}
Let us choose, for any $h\in\{1,\ldots,m\}$:
\[
\lambda= \sqrt{\frac{2N\log ({m}/{\varepsilon})}{ (1+ {1}/{k}
)^{2}\mathbf{P} [ (\chi^{\theta_{h}}
- \mathbf{P} \chi^{\theta_{h}} )^{2} ]
}} ,
\]
the bound becomes:
\begin{eqnarray*}
\lefteqn{\frac{1}{kN}\sum_{i=N+1}^{(k+1)N} \theta_{h}(X_{i})Y_{i}
-\frac{1}{N}\sum_{i=1}^{N} \theta_{h}(X_{i})Y_{i}}
\\
\lefteqn{\quad \leq \biggl(1+\frac{1}{k} \biggr) \biggl[ 2 \sqrt{\frac{ \mathbf{P}
[ (\chi^{\theta_{h}} - \mathbf{P}
\chi^{\theta_{h}} )^{2} ] \log ({m}/{\varepsilon})
}{ 2N }}
+ \frac{\mathcal{S}_{id}(\theta_{h})^{3}\log ({m}/{\varepsilon})} {
3N \mathbf{P} [ (\chi^{\theta_{h}} - \mathbf{P}
\chi^{\theta_{h}} )^{2} ]} \biggr] .}
\end{eqnarray*}

We use the reverse inequality exactly in the same way, we then
combine both inequality by a union bound argument and obtain the
following result. For any $\varepsilon>0$, with $P_{(k+1)N}$
probability at least $1-\varepsilon$ we have, for any
$h\in\{1,\ldots,m\}$:
\begin{eqnarray}\label{ancientheorem}
\nonumber r_{2} \bigl(\mathcal{C}^{h}\alpha_{1}^{h}\theta_{h} \bigr) -
r_{2} \bigl(\alpha_{2}^{h}\theta_{h} \bigr) &\leq&
\frac{ (1+ {1}/{k} )^{2}}{(1/(kN))\sum_{i=N+1}^{(k+1)N}
\theta_{h}(X_{i})^{2}} \biggl[
\frac{2\mathdsV _{\theta_{h}}\log (2m/{\varepsilon})}{N}
\\
&&{} + \frac{2 (\log ({2m}/{\varepsilon}) )^{{3}/{2}}
\mathcal{S}_{id}(\theta_{h})^{3}}{3 N^{{3}/{2}}
\mathdsV _{\theta_{h}}^{{1}/{2}} }
+ \frac{ (\log (2m/\varepsilon) )^{2}
\mathcal{S}_{id}(\theta_{h})^{6}}{ 9 N^{2}
\mathdsV _{\theta_{h}}^{2} }\biggr] ,
\end{eqnarray}
remember that:
\[
\mathdsV _{\theta} =
\mathbf{P} \bigl\{ \bigl[ \bigl(\theta(X)Y \bigr)
- \mathbf{P} \bigl(\theta(X)Y \bigr) \bigr]^{2} \bigr\}.
\]

We now give a new lemma.

\begin{lemma}
Let us assume that $P$ is such that, for any $h\in\{1,\ldots,m\}$:
\[
\exists\beta_{h} > 0,\ \exists B_{h} \geq0,\quad P \exp\bigl(\beta_{h} \big|\theta
_{h}(X)Y \big| \bigr) \leq B_{h} .
\]
This is for example the case if $\theta_{h}(X_{i})Y_{i}$ is
subgaussian, with any $\beta_{h} > 0$ and
\[
B_{h} = 2\exp\biggl\{\frac{\beta_{h}^{2}}{2} P \bigl[ \bigl(\theta_{h}(X)Y \bigr)^{2} \bigr] \biggr\}.
\]
Then we have, for any $\varepsilon\geq0$:
\[
P_{(k+1)N} \biggl\{\sup_{1\leq i\leq(k+1)N} \theta_{h}(X_{i})Y_{i}
\leq\frac{1}{\beta_{h}}\log\frac{(k+1)N B_{h}}{\varepsilon}
\biggr\} \geq1-\varepsilon.
\]
\end{lemma}

\begin{pf}
We have:
\begin{eqnarray*}
P_{(k+1)N} \Bigl(\sup_{1\leq i\leq(k+1)N} \theta_{h}(X_{i})Y_{i}
\geq s \Bigr)
&=&  P_{(k+1)N} \bigl(\exists i\in\bigl\{1,\ldots,(k+1)N\bigr\},
\theta_{h}(X_{i})Y_{i} \geq s \bigr)
\\
&=&  \sum_{i=1}^{(k+1)N} P 1_{\theta_{h}(X_{i})Y_{i} \geq s}
\\
&\leq& (k+1)N P \exp\bigl(\beta_{h} \big|\theta_{h}(X_{i})Y_{i} -
s \big| \bigr) \leq(k+1)N B_{h} \exp(-\beta_{h} s ) .
\end{eqnarray*}
Now, let use choose:
\[
s=\frac{1}{\beta_{h}} \log\frac{(k+1)N B_{h}}{\varepsilon} ,
\]
and we obtain the lemma.
\end{pf}

As a consequence, using a union bound argument, we have, for any
$\varepsilon\geq0$, with probability at least $1-\varepsilon$, for
any $h\in\{1,\ldots,m\}$:
\begin{eqnarray*}
\mathcal{S}_{id}(\theta_{h}) &=&  \sup_{i\in\{1,\ldots,(k+1)N\}}
\theta_{h}(X_{i})Y_{i} -\inf_{i\in\{1,\ldots,(k+1)N\}}
\theta_{h}(X_{i})Y_{i}
\leq \frac{2}{\beta_{h}} \log\frac{2 (k+1) m N B_{h}}{\varepsilon} .
\end{eqnarray*}

By plugging the lemma into Eq.~(\ref{ancientheorem}) we obtain
the theorem.
\end{pf*}

\subsection[Proof of Theorem 2.1: integration of the
transductive results]{Proof of Theorem \textup{\protect\ref{lastTH}}: integration of the
transductive results}

\label{subprooflastTH}

Actually, the proof is quite direct now: instead of using the
techniques given in the section devoted to the inductive case, we
use a result valid in the transductive case and integrate it with
respect to the test sample. This idea is quite classical in learning
theory, and was actually one of the reason for the introduction of
the transductive setting (see \cite{Vapnik} for example).
There are several ways to perform this integration (see for example
\cite{Classif}), here we choose to apply a result obtained by
Panchenko \cite{Panchenko} that gives a particularly simple result
here.

\begin{lemma}[(\cite{Panchenko}, Corollary~1)]
Let us assume that we have i.i.d. variables $T_{1},\ldots,T_{N}$ (with
distribution $P$ and values in $\mathdsR $) and an independent copy
$T'=(T'_{1},\ldots,T'_{N})$ of $T=(T_{1},\ldots,T_{N})$. Let
$\xi_{j}(T,T')$ for $j\in\{1,2,3\}$ be three measurables functions
taking values in $\mathdsR $, and $\xi_{3}\geq0$. Let us assume
that we know two constants $A\geq1$ and $a>0$ such that, for any
$u>0$:
\[
P^{\otimes2N} \bigl[\xi_{1}\bigl(T,T'\bigr)\geq\xi_{2}\bigl(T,T'\bigr)+\sqrt{\xi_{3}\bigl(T,T'\bigr)u} \bigr]
\leq A\exp(-au) .
\]
Then, for any $u>0$:
\[
P^{\otimes2N} \bigl\{P^{\otimes2N} \bigl[\xi_{1}\bigl(T,T'\bigr)|T \bigr]
\geq P^{\otimes2N} \bigl[\xi_{2}\bigl(T,T'\bigr)|T \bigr]+
\sqrt{P^{\otimes2N} \bigl[\xi_{3}\bigl(T,T'\bigr)|T \bigr]u} \bigr\} \leq A\exp(1-au) .
\]
\end{lemma}

\begin{pf*}{Proof of Theorem \ref{lastTH}}
A simple application of the first inequality of Lemma
\ref{learninglemma} (given as a tool for the proof of the
transductive results) with
$\varepsilon>0$, any $k\in\{1,\ldots,m\}$, $g=id$,
$\eta=1+\log\frac{2m}{\varepsilon}$ and:
\[
\lambda_{k} = \sqrt{\frac{N\eta}{(1/N)\sum_{i=1}^{2N}\theta
_{k}(X_{i})^{2}Y_{i}^{2}}}
\]
leads us to the following bound, for any $k$:
\[
P^{\otimes2N} \exp\biggl[\sqrt{N\eta}
\frac{({1}/{N})\sum_{i=1}^{N} [\theta_{k}(X_{i})Y_{i}-\theta
_{k}(X_{i+N})Y_{i+N} ]}
{\sqrt{({1}/{N})\sum_{i=1}^{2N}\theta_{k}(X_{i})^{2}Y_{i}^{2}}}-2\eta\biggr]
\leq\exp(-\eta),
\]
or:
\[
P^{\otimes2N} \Biggl[
\frac{1}{N}\sum_{i=1}^{N} \bigl[\theta_{k}(X_{i})Y_{i}-\theta
_{k}(X_{i+N})Y_{i+N} \bigr]
\geq
\sqrt{\frac{4\eta}{N^{2}}\sum_{i=1}^{2N}\theta_{k}(X_{i})^{2}Y_{i}^{2}}
\Biggr]
\leq\exp(-\eta) = \frac{\varepsilon}{2k\exp(1)}.
\]
We now apply Panchenko's lemma with:
\begin{eqnarray*}
\lefteqn{T_{i} = \theta_{k}(X_{i})Y_{i},\qquad
T_{i}'=\theta_{k}(X_{i+N})Y_{i+N},}
\\
\lefteqn{\xi_{1}\bigl(T,T'\bigr) = \frac{1}{N}\sum_{i=1}^{N}T_{i} ,\qquad  \xi_{2}\bigl(T,T'\bigr) =
\frac{1}{N}\sum_{i=1}^{N}T_{i}',}
\\
\lefteqn{\xi_{3}\bigl(T,T'\bigr) = \frac{2}{N^{2}}\sum_{i=1}^{2N}\theta
_{k}(X_{i})^{2}Y_{i}^{2} \geq0,}
\end{eqnarray*}
and $A=a=1$. We obtain:
\begin{eqnarray*}
P^{\otimes2N} \Biggl[
\frac{1}{N}\sum_{i=1}^{N} \bigl[\theta_{k}(X_{i})Y_{i}-P \bigl[\theta_{k}(X)Y \bigr] \bigr]
\geq
\sqrt{\frac{4\eta}{N^{2}}\sum_{i=1}^{N} \bigl[\theta
_{k}(X_{i})^{2}Y_{i}^{2}+P \bigl[\theta_{k}(X)^{2}Y^{2} \bigr] \bigr]}
\Biggr]
\leq\exp(1-\eta) = \frac{\varepsilon}{2k}.
\end{eqnarray*}
Remark finally that:
\[
P \bigl[\theta_{k}(X)^{2}Y^{2} \bigr] \leq P \bigl[\theta_{k}(X)^{2} \bigr]\bigl(B^{2}+\sigma^{2}\bigr).
\]
We proceed exactly in the same way with the reverse inequalities for
any $k$ and combine the obtained $2m$ inequalities to obtain the
result:
\begin{eqnarray*}
\lefteqn{P^{\otimes N} \Biggl\{ \exists k\in\{1,\ldots,m\},
\frac{1}{N}\sum_{i=1}^{N} \big|\theta_{k}(X_{i})Y_{i}-P \bigl[\theta_{k}(X)Y \bigr]
\big|}
\\
\lefteqn{\quad \geq
\sqrt{\frac{4+4\log ({2m}/{\varepsilon})}{N^{2}}\sum_{i=1}^{N} \bigl\{\theta
_{k}(X_{i})^{2}Y_{i}^{2}+P \bigl[\theta_{k}(X)^{2} \bigr]\bigl(B^{2}+\sigma^{2}\bigr) \bigr\}}
\Biggr\}}
\\
\lefteqn{\qquad =P^{\otimes2N} \Biggl\{ \exists k\in\{1,\ldots,m\},
\frac{1}{N}\sum_{i=1}^{N} \big|\theta_{k}(X_{i})Y_{i}-P \bigl[\theta_{k}(X)Y \bigr]
\big|}
\\
\lefteqn{\quad \geq
\sqrt{\frac{4+4\log ({2m}/{\varepsilon})}{N^{2}}\sum_{i=1}^{N} \bigl\{\theta
_{k}(X_{i})^{2}Y_{i}^{2}+P \bigl[\theta_{k}(X)^{2} \bigr]\bigl(B^{2}+\sigma^{2}\bigr) \bigr\}}
\Biggr\} \leq\varepsilon}
\end{eqnarray*}
that ends the proof.
\end{pf*}

\subsection[Proof of Theorems 4.1 and 4.2: Theorem 2.3 used as an oracle inequality]{Proof of
Theorems \textup{\protect\ref{speed}} and \textup{\protect\ref{rate2}}: Theorem \textup{\protect\ref{propalgo}}
used as an oracle inequality}
\label{subproofrate}

\begin{pf*}{Proof of Theorem \ref{speed}}
Let us begin the proof with a general $m$ and $\varepsilon$, the
reason of the choice $m=N$ and $\varepsilon=N^{-2}$ will become
clear. Let us also call $\mathcal{E}(\varepsilon)$ the event
satisfied with probability at least $1-\varepsilon$ in Theorem
\ref{lastTH}. We have:
\[
P^{\otimes
N} \bigl[ \big\|\varPi_{P}^{\mathcal{F},m}\hat{\theta}-f \big\|_{P}^{2} \bigr]
= P^{\otimes
N} \bigl[1_{\mathcal{E}(\varepsilon)} \big\|\varPi_{P}^{\mathcal{F},m}\hat
{\theta}-f \big\|_{P}^{2} \bigr]
+ P^{\otimes
N} \bigl[ (1-1_{\mathcal{E}(\varepsilon)} ) \big\|\varPi_{P}^{\mathcal
{F},m}\hat{\theta}-f \big\|_{P}^{2} \bigr].
\]
First of all, it is obvious that:
\begin{eqnarray*}
P^{\otimes
N} \bigl[ (1-1_{\mathcal{E}(\varepsilon)} ) \big\|\varPi_{P}^{\mathcal
{F},m}\hat{\theta}-f \big\|_{P}^{2} \bigr]
&\leq& 2 P^{\otimes
N} \bigl[ (1-1_{\mathcal{E}(\varepsilon)} )
\bigl( \big\|\varPi_{P}^{\mathcal{F},m}\hat{\theta} \big\|_{P}^{2}+ \|f \|_{P}^{2} \bigr) \bigr]
\\
&\leq& 2 \varepsilon \bigl(B^{2}m + B^{2} \bigr) =
2\varepsilon(m+1)B^{2}.
\end{eqnarray*}
For the other term, just remark that, for any $m'\leq m$:
\begin{eqnarray*}
\big\|\varPi_{P}^{\mathcal{F},N}\hat{\theta}-f \big\|_{P}^{2} &=&
\big\|\varPi_{P}^{\mathcal{F},m}\varPi_{P}^{m,\varepsilon}\cdots\varPi
_{P}^{1,\varepsilon}0-f \big\|_{P}^{2}
\leq\big\|\varPi_{P}^{m,\varepsilon}\cdots\varPi_{P}^{1,\varepsilon}0-f \big\|_{P}^{2}
\leq
\big\|\varPi_{P}^{m',\varepsilon}\cdots\varPi_{P}^{1,\varepsilon}0-f \big\|_{P}^{2}
\\
&\leq& \sum_{k=1}^{m'}
\frac{4 [1+\log ({2m}/{\varepsilon}) ]}{N} \Biggl[
\frac{1}{N}\sum_{i=1}^{N}\theta_{k}(X_{i})^{2}Y_{i}^{2} + B^{2} +
\sigma^{2} \Biggr] +
\|\overline{\theta}_{m'}-f \|_{P}^{2}.
\end{eqnarray*}
This is where Theorem \ref{propalgo} has been used as an oracle
inequality: the estimator that we have, with $m\geq m'$, is better
than the one with the ``good choice'' $m'$. We also have:
\begin{eqnarray*}
P^{\otimes
N} \bigl[1_{\mathcal{E}(\varepsilon)} \big\|\varPi_{P}^{\mathcal{F},m}\hat
{\theta}-f \big\|_{P}^{2} \bigr]
&\leq& P^{\otimes N} \Biggl[ \sum_{k=1}^{m'}
\frac{4 [1+\log ({2m}/{\varepsilon})]}{N} \Biggl[
\frac{1}{N}\sum_{i=1}^{N}\theta_{k}(X_{i})^{2}Y_{i}^{2} + B^{2} +
\sigma^{2} \Biggr] \Biggr] + \bigl(m'\bigr)^{-2\beta}C
\\
&\leq& m'\frac{8 [1+\log (2m/\varepsilon) ]}{N}
\bigl[B^{2} + \sigma^{2} \bigr].
\end{eqnarray*}
So finally, we obtain, for any $m'\leq m$:
\[
P^{\otimes
N} \bigl[ \big\|\varPi_{P}^{\mathcal{F},m}\hat{\theta}-f \big\|_{P}^{2} \bigr]
\leq m' \frac{8 [1+\log (2m/\varepsilon) ]}{N}
\bigl[B^{2} + \sigma^{2} \bigr]
+ \bigl(m'\bigr)^{-2\beta}C +2\varepsilon(m+1)B^{2}.
\]
The choice of:
\[
m'= \biggl(\frac{N}{\log N} \biggr)^{{1}/{(2\beta+1)}}
\]
leads to a first term of order $N^{{-2\beta}/{(2\beta+1)}}\log
\frac{m}{\varepsilon}(\log N)^{{2\beta}/{(2\beta+1)}}$ and a
second term of order $N^{{-2\beta}/{(2\beta+1)}}\times (\log
N)^{{2\beta}/{(2\beta+1)}}$. The choice of $m=N$ and
$\varepsilon=N^{-2}$ gives a first and a second term of the desired
order $N^{{-2\beta}/{(2\beta+1)}}(\log
N)^{{2\beta}/{(2\beta+1)}}$ while keeping the third term at order
$N^{-1}$. This proves the theorem.
\end{pf*}

\begin{pf*}{Proof of Theorem \ref{rate2}}
Here again let us write $\mathcal{E}(\varepsilon)$ the event
satisfied with probability at least $1-\varepsilon$ in Theorem~\ref{lastTH}. We have:
\[
P^{\otimes N}
\bigl[ \big\|\varPi_{P}^{\mathcal{F},m}\hat{\theta}-f \big\|_{P}^{2} \bigr]
= P^{\otimes N} \bigl[1_{\mathcal{E}(\varepsilon)}
\big\|\varPi_{P}^{\mathcal{F},m}\hat{\theta}-f \big\|_{P}^{2} \bigr]
+ P^{\otimes
N} \bigl[ (1-1_{\mathcal{E}(\varepsilon)} )
\big\|\varPi_{P}^{\mathcal{F},m}\hat{\theta}-f \big\|_{P}^{2}
\bigr].
\]
For the first term we still have:
\[
\big\|\varPi_{P}^{\mathcal{F},m}\hat{\theta}-f \big\|_{P}^{2} \leq
2(m+1)B^{2}.
\]
For the second term, let us write the expansion of $f$ into our
wavelet basis:
\[
f=\alpha\phi+ \sum_{j=0}^{\infty}\sum_{k=1}^{2^{j}}\beta_{j,k}\psi
_{j,k} ,
\]
and
\[
\hat{\theta}(x) = \tilde{\alpha} \phi+ \sum_{j=0}^{J}\sum
_{k=1}^{2^{j}}\tilde{\beta}_{j,k}\psi_{j,k}
\]
the estimator $\hat{\theta}$. Let us put
$ J=2^{\lfloor (\log N)/{\log2} \rfloor}$.
\begin{eqnarray*}
\big\|\varPi_{P}^{\mathcal{F},m}\hat{\theta}-f \big\|_{P}^{2} &\leq&
\|\hat{\theta}-f \|_{P}^{2} =
\big\|\varPi_{P}^{m,\varepsilon}\cdots\varPi_{P}^{1,\varepsilon}0-f \big\|_{P}^{2}
\\
&=&  (\tilde{\alpha} - \alpha)^{2} +
\sum_{j=0}^{J}\sum_{k=1}^{2^{j}}(\tilde{\beta}_{j,k}-\beta_{j,k})^{2}
+ \sum_{j=J+1}^{\infty}\sum_{k=1}^{2^{j}}\beta_{j,k}^{2}
\\
&\leq& (\tilde{\alpha} - \alpha)^{2} +
\sum_{j=0}^{J}\sum_{k=1}^{2^{j}}(\tilde{\beta}_{j,k}-\beta
_{j,k})^{2}1\bigl(|\beta_{j,k}|\geq
\kappa\bigr) + \sum_{j=0}^{J}\sum_{k=1}^{2^{j}} \beta_{j,k}^{2}
1\bigl(|\beta_{j,k}|< \kappa\bigr)
+ \sum_{j=J+1}^{\infty}\sum_{k=1}^{2^{j}}\beta_{j,k}^{2}
\end{eqnarray*}
for any $\kappa\geq0$, as soon as $\mathcal{E}(\varepsilon)$ is
satisfied (here again we used Theorem \ref{propalgo} as an oracle
inequality). Now, we follow the technique used in \cite{Ondel2} and
\cite{Ondel} (see also the end of the third chapter in \cite{Cat7}).
As soon as $\mathcal{E}(\varepsilon)$ is satisfied we have:
\begin{eqnarray*}
\sum_{j=0}^{J}\sum_{k=1}^{2^{j}}(\tilde{\beta}_{j,k}-\beta
_{j,k})^{2}1\bigl(|\beta_{j,k}| \geq
\kappa\bigr) &\leq& \frac{8(B^{2}+\sigma^{2})
\log(2m/{\varepsilon})}{N}
\sum_{j=0}^{J}\sum_{k=1}^{2^{j}}1\bigl(|\beta_{j,k}|\geq\kappa\bigr)
\\
&\leq& \frac{8(B^{2}+\sigma^{2}) \log (2m/\varepsilon)}{N}
\sum_{j=0}^{J}\sum_{k=1}^{2^{j}}
\biggl(\frac{|\beta_{j,k}|}{\kappa} \biggr)^{{2}/{(2s+1)}}
\\
&=&  \frac{8(B^{2}+\sigma^{2}) \log(2m/\varepsilon)}{N}
\kappa^{-{2}/{(2s+1)}} \sum_{j=0}^{J}\sum_{k=1}^{2^{j}}
|\beta_{j,k}|^{{2}/{(2s+1)}}.
\end{eqnarray*}
In the same way, we have:
\[
\sum_{j=0}^{J}\sum_{k=1}^{2^{j}} \beta_{j,k}^{2}
1\bigl(|\beta_{j,k}|< \kappa\bigr) \leq\kappa^{2-{2}/{(1+2s)}}
\sum_{j=0}^{J}\sum_{k=1}^{2^{j}}
|\beta_{j,k} |^{{2}/{(1+2s)}} .
\]
So we have to give an upper bound on the quantity:
\[
\sum_{j=0}^{J}\sum_{k=1}^{2^{j}} |\beta_{j,k}|^{{2}/{(2s+1)}} .
\]
By H\"{o}lder's inequality we have, as soon as $p\geq\frac{2}{2s+1}$:
\[
\sum_{j=0}^{J}\sum_{k=1}^{2^{j}} |\beta_{j,k}|^{{2}/{(2s+1)}} \leq
\sum_{j=0}^{J} \Biggl[2^{j (1+ {1}/{2}- {1}/{p} )}
\sum_{k=1}^{2^{j}} |\beta_{j,k}|^{p} \Biggr]^{{2}/{(1+2s)}} \leq
\|f\|_{s,p,q}^{{2}/{(1+2s)}}
J^{ (1-{2}/{((1+2s)q)} )_{+}} ,
\]
let us put
$C'=\|f\|_{s,p,q}^{{2}/{(1+2s)}}$. Finally, note that we have, for
$p\geq2$:
\[
\sum_{j=J+1}^{\infty}\sum_{k=1}^{2^{j}} \beta_{j,k}^{2} \leq
\sum_{j=J+1}^{\infty} \Biggl(\sum_{k=1}^{2^{j}}
\beta_{j,k}^{p} \Biggr)^{{2}/{p}}
2^{j (1-{2}/{p} )}.
\]
As $f\in B_{s,p,q} \subset B_{s,p,\infty}$ we have:
\[
\Biggl(\sum_{k=1}^{2^{j}} \beta_{j,k}^{p} \Biggr)^{{2}/{p}} \leq C' 2^{-2j
(s+ {1}/{2}- {1}/{p} )}
\]
for some $C''$ and so:
\[
\sum_{j=J+1}^{\infty}\sum_{k=1}^{2^{j}} \beta_{j,k}^{2} \leq C'''
2^{-2Js}
\]
for some $C'''$. In the case where $p<2$ we use (see \cite{Ondel},
for $s>\frac{1}{p}-\frac{1}{2}$):
\[
B_{s,p,q} \subset B_{s-{1}/{p}+{1}/{2},2,q}
\]
to obtain:
\[
\sum_{j=J+1}^{\infty}\sum_{k=1}^{2^{j}} \beta_{j,k}^{2} \leq
C''''2^{-2J (s+ {1}/{2}- {1}/{p} )}
\leq C'''' 2^{-J} .
\]
So we have:
\begin{eqnarray*}
P^{\otimes N} d^{2}(\tilde{f},f) &\leq& 2(m+1)\varepsilon
\bigl(B^{2}+\sigma^{2}\bigr)
+ \frac{8(B^{2}+\sigma^{2})\log (2m/{\varepsilon})}{N} \bigl(1+C'\kappa
^{-{2}/{(1+2s)}}J^{ (1-{2}/{((1+2s)q)} )_{+}} \bigr)
\\
&&{} +C'\kappa^{2-{2}/{(1+2s)}}
J^{ (1-{2}/{((1+2s)q)} )_{+}}
+C''' \bigl(2^{-J} \bigr)^{2s} + C''''2^{-J} .
\end{eqnarray*}
Let us remember that:
\[
\frac{N}{2} \leq m = 2^{J} \leq N
\]
and that $\varepsilon=N^{-2}$, and take:
\[
\kappa=\sqrt{\frac{\log N}{N}}
\]
to obtain the desired rate of convergence.
\end{pf*}

\section*{Acknowledgments}

I would like to thank my PhD advisor, Professor Olivier Catoni, for his
constant help, and the anonymous referee for
very useful comments and remarks.

\end{document}